\documentclass[pdflatex,sn-mathphys-num]{sn-jnl}%

\usepackage{graphicx}%
\usepackage{multirow}%
\usepackage{amsmath,amssymb,amsfonts}%
\usepackage{amsthm}%
\usepackage{mathrsfs}%
\usepackage[title]{appendix}%
\usepackage{xcolor}%
\usepackage{textcomp}%
\usepackage{manyfoot}%
\usepackage{booktabs}%
\usepackage{tikz}
\usepackage{pifont} % for \ding command
\usepackage{algorithm}%
\usepackage{algorithmicx}%
\usepackage{algpseudocode}%
\usepackage{listings}%

\usepackage{float}
\usepackage{subfig}
\usepackage{bm}
\usepackage{caption} 
\usepackage{subcaption}
\usepackage[capitalize,noabbrev]{cleveref}
\usepackage{tabularx}
\theoremstyle{thmstyleone}%
\theoremstyle{thmstyletwo}%

\theoremstyle{thmstylethree}%
\newcommand{\mat}[1]{{\bf #1}}

\newcommand{\ba}{\mathbf a}
\newcommand{\bb}{\mathbf b}

\newcommand{\bh}{\mathbf h}

\newcommand{\bw}{\mathbf w}
\newcommand{\bx}{\mathbf x}

\newcommand{\bz}{\mathbf z}

\newcommand{\calb}{\mathcal B}

\newcommand{\caln}{\mathcal N}

\newcommand{\R}{\mathbb R}

\newcommand{\dho}{\beta}

\newcommand{\rev}[1]{\textcolor{black}{#1}}

\begin{document}

\title[ELM-FBPINN]{ELM-FBPINNs: An Efficient Multilevel Random Feature Method}

\author[1]{\fnm{Samuel} \sur{Anderson}}\email{sam.anderson@strath.ac.uk}

\author*[2]{\fnm{Victorita} \sur{Dolean}}\email{v.dolean.maini@tue.nl}
\equalcont{These authors contributed equally to this work.}

\author[3]{\fnm{Ben} \sur{Moseley}}\email{b.moseley@imperial.ac.uk}
\equalcont{These authors contributed equally to this work.}

\author[1]{\fnm{Jennifer} \sur{Pestana}}\email{jennifer.pestana@strath.ac.uk}
\equalcont{These authors contributed equally to this work.}

\affil*[1]{\orgdiv{Department of Mathematics and Statistics}, \orgname{University of Strathclyde}, \orgaddress{\city{Glasgow}, \postcode{G1 1XH}, \country{UK}}}

\affil[2]{\orgdiv{Department of Mathematics and Computer Science}, \orgname{Eindhoven University of Technology}, \orgaddress{\city{Eindhoven}, \postcode{5600 MB}, \country{The Netherlands}}}

\affil[3]{\orgdiv{Department of Earth Science and Engineering}, \orgname{Imperial College London}, \orgaddress{\city{London}, \postcode{SW7 2AZ},  \country{UK}}}


\abstract{
Domain-decomposed variants of physics-informed neural networks (PINNs) such as finite basis PINNs (FBPINNs) mitigate some of PINNs' issues like slow convergence and spectral bias through localisation, but still rely on iterative nonlinear optimisation within each subdomain. In this work, we propose a hybrid approach that combines multilevel domain decomposition and partition-of-unity constructions with random feature models, yielding a method referred to as multilevel ELM-FBPINN. By replacing trainable subdomain networks with extreme learning machines, the resulting formulation eliminates backpropagation entirely and reduces training to a structured linear least-squares problem. 
We provide a systematic numerical study comparing ELM-FBPINNs and multilevel ELM-FBPINNs with standard PINNs and FBPINNs on representative benchmark problems, demonstrating that ELM-FBPINNs and multilevel ELM-FBPINNs achieve competitive \rev{errors} while significantly accelerating convergence and improving robustness with respect to architectural and optimisation parameters. Through ablation studies, we further clarify the distinct roles of domain decomposition and random feature enrichment in controlling expressivity, conditioning, and scalability.
}




\maketitle

\section{Introduction}
\label{sec:introduction}

Accurate and efficient numerical solution of partial differential equations (PDEs) remains a central challenge in computational science and engineering. Classical discretisation techniques such as finite difference and finite element methods have matured into highly reliable tools, yet they can become costly or difficult to apply in problems involving high dimensionality, complex geometries, or strongly multiscale behaviour. These challenges have motivated increasing interest in alternative computational paradigms capable of combining flexible function approximation with physical constraints. In recent years, \emph{scientific machine learning} (SciML) has emerged as a promising framework for tackling such problems. Among the various approaches proposed in this area, \emph{physics-informed neural networks} (PINNs) \cite{Lagaris1998,Raissi2019} have attracted considerable attention. In PINNs, the solution of a PDE is approximated by a neural network whose parameters are trained by minimising a loss function that penalises the residual of the governing PDE and its boundary conditions. This formulation allows PDE solving to be recast as a continuous optimisation problem and offers several appealing features, including mesh-free discretisation, straightforward implementation through automatic differentiation, and a unified framework for forward and inverse problems.

Despite these advantages, standard PINNs trained via gradient-based optimisation exhibit several well-documented difficulties. Training can be slow and sensitive to hyperparameters, and the resulting optimisation problem is often poorly conditioned. Moreover, neural networks trained by gradient descent are known to display a \emph{spectral bias} toward low-frequency functions \cite{Rahaman2018}, which limits their ability to resolve highly oscillatory or multiscale solutions. As a result, classical PINNs frequently struggle with problems containing high-frequency content or widely separated spatial scales \cite{Moseley2023}. To address these limitations, a variety of domain decomposition strategies have been proposed. Conservative PINNs (cPINNs) \cite{jagtap2020conservative} enforce flux continuity across non-overlapping subdomains, while extended PINNs (XPINNs) \cite{jagtap2020extended} allow flexible space–time decompositions. 

A particularly effective approach is the \emph{Finite Basis Physics-Informed Neural Network} (FBPINN) method introduced by Moseley et al.\ \cite{Moseley2023}. In FBPINNs, the computational domain is decomposed into overlapping subdomains, each equipped with a small neural network whose output is multiplied by a smooth window function forming a partition of unity. The global solution is obtained by summing the local contributions across all subdomains. This localisation substantially mitigates spectral bias by reducing the complexity of the function each network must approximate \cite{dolean2022finite, Dolean:MDD:2024}. Nevertheless, FBPINNs still rely on gradient-based training of all network parameters, which remains computationally expensive and can limit robustness.

An alternative strategy is to replace fully trainable neural networks with \emph{random feature models}. In particular, \emph{Extreme Learning Machines} (ELMs) \cite{Huang2006} fix the hidden-layer weights at random and determine only the output weights by solving a linear least-squares problem. This idea transforms training from a nonlinear optimisation problem into a linear algebra problem and is closely related to the \emph{Random Feature Method} (RFM) \cite{chen:2022:BTM}. Randomised neural network models have a long history, including random vector functional-link networks \cite{pao1992functional} and reservoir computing \cite{Jaeger:2009:RCA}, and have recently been applied to PDE solving. Examples include the Local Extreme Learning Machine method of Dong and Li \cite{DongLi2021_ELM-DDM}, distributed learning machine approaches \cite{Dwivedi2021}, and the Extreme Theory of Functional Connections framework \cite{Schiassi2021_ETFC}. These methods demonstrate that fixed random features can provide competitive accuracy while dramatically simplifying training \cite{chen:2022:BTM, Anderson:2024:ELM, datar2026frozenpinn}.

In this work, we investigate the combination of these two ideas: domain localisation and random feature approximation. Specifically, we replace the trainable subdomain networks in the FBPINN framework with extreme learning machines, yielding a method we refer to as \emph{ELM-FBPINN}. A preliminary study was performed in the first version of \cite{Anderson:2024:ELM}. In the current work, we go one step further by pushing deeper the study of its performance and exploring numerically both a one-level method and its multilevel generalisation. The resulting formulation retains the partition-of-unity domain decomposition structure of FBPINNs while eliminating gradient-based optimisation entirely. Instead, the PDE solution is obtained by solving a global linear least-squares system for the output weights associated with randomly initialised basis functions.

This reformulation raises several fundamental questions regarding the role of optimisation, localisation, and model capacity in physics-informed learning. In particular, we seek to understand whether replacing nonlinear optimisation with linear least-squares solves can improve robustness and convergence behaviour without sacrificing accuracy, and how the resulting method behaves as problem difficulty and model capacity vary. The present study is therefore guided by the following research questions:
\begin{itemize}
\item \textbf{RQ1:} Does replacing gradient-based optimisation with a linear least-squares formulation improve convergence speed and robustness while maintaining solution accuracy?
\item \textbf{RQ2:} To what extent does domain localisation mitigate spectral bias and enable accurate approximation of highly oscillatory or multiscale solutions?
\item \textbf{RQ3:} How does the number of random basis functions per subdomain influence model expressivity, redundancy, and the conditioning of the resulting least-squares system?
\item \textbf{RQ4:} Can the method maintain stable performance as problem difficulty increases when model capacity is scaled accordingly, and can multilevel modeling improve this performance?
\end{itemize}

To address these questions, we compare three architectures: a standard PINN, an FBPINN, and the proposed ELM-FBPINN. This comparison allows us to isolate the effects of domain decomposition and training strategy independently. The methods are evaluated on \rev{three} benchmark problems: a one-dimensional damped harmonic oscillator, which exhibits highly oscillatory temporal dynamics; a two-dimensional multi-scale Laplacian problem designed to test spatial multiscale behaviour; \rev{and a two-dimensional inhomogeneous Helmholtz problem on a non-convex L-shaped domain with oscillatory, non-separable spatial structure.}

Through these experiments, we aim to clarify the respective roles of localisation, random feature representations, and linear least-squares formulations in physics-informed neural networks. More broadly, the results contribute to an emerging perspective in which neural PDE solvers are viewed not only as trainable black-box models but also as approximation schemes whose design choices directly influence their numerical behaviour.

\section{Related work}
The present work can be viewed as a generalisation of the random feature formulation introduced in the first version of  \cite{Anderson:2024:ELM}, extending it to a multilevel domain-decomposition setting. It is also closely related to the Random Feature Method (RFM) proposed in \cite{chen:2022:BTM}. The main differences with respect to \cite{chen:2022:BTM} lie in the treatment of domain localisation and the exploitation of sparsity. In the present work, the neural basis functions are organised through a domain-decomposition framework and evaluated only on the collocation points associated with their corresponding subdomains. This leads to a naturally sparse system matrix, which we explicitly exploit by assembling the matrix in sparse form and solving the resulting least-squares problem using sparse iterative solvers. In contrast, existing RFM formulations typically construct dense systems and do not exploit the sparsity induced by localisation.

Several recent works have also explored random-feature and frozen-network approaches for PDE solvers, including multiscale and domain-decomposition settings. In particular, \cite{yildiz2026fastmultiscale} proposes a multiscale PDE solver based on domain decomposition and random features, and \cite{datar2026frozenpinn} introduces frozen physics-informed neural networks to accelerate training through linear least-squares formulations. Finally, other works have investigated preconditioning strategies for random-feature-based PDE solvers. In particular, \cite{Shang_Heinlein_Mishra_Wang_2025} and \cite{vanBeek2026RRQR} propose preconditioning techniques aimed at improving the conditioning of the resulting least-squares systems and do not consider multilevel formulations. 

While these works share the same general philosophy of replacing gradient-based training with linear solvers, the present study provides a more thorough assessment of the multilevel domain-decomposition framework, with a particular focus on approximation properties, sparsity structure, conditioning, and scalability. 
In particular, we explicitly study how localisation, multilevel decomposition, and basis enrichment interact, and how these design choices affect accuracy, robustness, and computational efficiency.

\section{The multilevel ELM-FBPINN method}

We now extend the random-feature FBPINN formulation introduced in \cite{Anderson:2024:ELM} and used in \cite{vanBeek2026RRQR} to a multilevel domain-decomposition setting. The goal is to retain the localization advantages of the FBPINN framework while replacing gradient-based training by a linear least-squares solve using randomised neural features.

\paragraph{Problem and model definition}

We start with a PDE defined in a domain $\Omega \subset \R^d $
\[ \mathcal{N}[u(\bx)] = f(\bx), \ \ \bx \in \Omega, \]
under boundary conditions given by 
\[ \mathcal{B}_b[u(\bx)] = g_b(\bx), \ \ \bx \in \Gamma_b \subset \partial \Omega, \ \ b=1,\dots,B, \]
where $\caln$ and the $\calb_b$ are assumed to be linear operators and $f$ and $g_b$ are forcing and boundary functions. Instead of using a single decomposition into $J$ subdomains, we introduce a hierarchy of decompositions indexed by levels $\ell = 1,\dots,L$. At each level $\ell$, the domain is covered by $J_\ell$ overlapping subdomains
\(
\Omega = \bigcup_{j=1}^{J_\ell} \Omega_{j}^{(\ell)}.
\)
Associated with these subdomains are compactly supported window functions
\(
\{\omega_{j}^{(\ell)}(\bx)\}_{j=1}^{J_\ell},
\)
which form a partition of unity at each level,
\[
\sum_{j=1}^{J_\ell} \omega_{j}^{(\ell)}(\bx) = 1,
\qquad \bx \in \Omega .
\]
\rev{
The window functions are defined on finite subdomains and are taken to be zero outside their prescribed support. At the domain boundaries, subdomains are truncated accordingly, and the partition-of-unity normalisation ensures that the sum of all window functions remains equal to one. This construction is standard in overlapping domain decomposition methods and guarantees a consistent global representation.
}

Typically the decompositions become progressively finer as the level increases, so that level $\ell=1$ corresponds to a coarse representation of the solution and higher levels introduce increasingly localised corrections.
Using this hierarchy, we approximate the PDE solution as a sum of localised random feature expansions across all levels,
\begin{equation}
\label{eq:multilevel_solution}
\hat{u}(\bx,\ba) =
\frac{1}{L}\sum_{\ell=1}^{L}
\sum_{j=1}^{J_\ell}
\omega_{j}^{(\ell)}(\bx)
\sum_{k=1}^{K}
a_{jk}^{(\ell)} \,
\phi_{jk}^{(\ell)}(\bx,\theta_{jk}^{(\ell)}).
\end{equation}
Here $\phi_{jk}^{(\ell)}(\bx,\theta_{jk}^{(\ell)})$ are neural network basis functions, $K$ denotes the number of random features per subdomain and  $a_{jk}^{(\ell)}$ are trainable output weights. As in the single-level formulation, the neural features are taken to be shallow fully connected networks
\begin{equation}
\begin{aligned}
\bz_{j}^{(\ell,0)} &= \bx, \\
\bz_{j}^{(\ell,m)} &= \Xi\!\left(W_{j}^{(\ell,m)}\bz_{j}^{(\ell,m-1)} + \bb_{j}^{(\ell,m)}\right),
\quad m=1,\dots,h-1, \\
\phi_{jk}^{(\ell)}(\bx,\theta_{jk}^{(\ell)})
&=
\Xi\!\left((\bw_{jk}^{(\ell)})^{\top}\bz_{j}^{(\ell,h-1)}\right),
\end{aligned}
\end{equation}
where $\Xi$ is a nonlinear activation function applied elementwise.
The parameters
\[
\theta_{jk}^{(\ell)} = \{W_{j}^{(\ell,m)}, \bb_{j}^{(\ell,m)}\}_{m=1}^{h-1} \cup \{\bw_{jk}^{(\ell)}\}
\]
are randomly initialised and then kept fixed. Consequently, only the coefficients $a_{jk}^{(\ell)}$ remain as unknown parameters.
This formulation therefore produces a multilevel random-feature approximation in which each level contributes localised neural basis functions with support determined by the window functions $\omega_j^{(\ell)}$.

\paragraph{Remark on multilevel coupling.}
The effect of the multilevel construction can be understood by analogy with multigrid methods and hierarchical basis representations. In a single-level decomposition, information propagates between distant subdomains only indirectly through overlapping window functions, which can lead to slow global coordination when the solution contains long-range or low-frequency components.  By introducing multiple levels, the approximation space is enriched in a hierarchical manner: coarse levels provide basis functions with large spatial support that capture the global structure of the solution, while finer levels introduce increasingly localised corrections that resolve high-frequency features. The partition-of-unity construction ensures that these contributions are smoothly combined across levels, yielding a stable global approximation. This results in a basis in which low-frequency components are represented on coarse levels and high-frequency components on fine levels, analogous to the separation of error modes in multigrid. In this sense, the multilevel ELM-FBPINN facilitates “global subdomain communication” by enabling long-range interactions to be represented explicitly at coarse scales, rather than relying solely on local overlap between neighbouring subdomains.
However, in contrast to multigrid methods, multilevel ELM-FBPINN does not follow a classical residual-correction structure. Instead, it can be interpreted as a hierarchical basis construction, where each level contributes functions with different spatial support, as described above.

\paragraph{Loss function and least-squares formulation}

The unknown coefficients are determined by minimizing the physics-informed loss functional
\begin{equation}
\label{eq:multilevel_loss}
\mathcal{L}(\bm{a}) =
\frac{1}{N_I}\sum_{n=1}^{N_I}
\left(
\mathcal{N}[\hat{u}(\bx_n,\bm{a})] - f(\bx_n)
\right)^2
+
\sum_{b=1}^{B}
\frac{\lambda_b}{N_B^{(b)}}
\sum_{i=1}^{N_B^{(b)}}
\left(
\mathcal{B}_b[\hat{u}(\bx_i^{(b)},\bm{a})] - g_b(\bx_i^{(b)})
\right)^2 .
\end{equation}
Here $\{\bx_n\}_{n=1}^{N_I}$ are interior collocation points and $\{\bx_i^{(b)}\}_{i=1}^{N_B^{(b)}}$ are boundary collocation points associated with boundary operators $\mathcal{B}_b$. 
Since the operators $\mathcal{N}$ and $\mathcal{B}_b$ are linear, the loss functional remains quadratic in the coefficients $\bm{a}$ even in the multilevel setting.
 Consequently, the training problem reduces to a linear least-squares problem. 
 
Let 
\(
N_\Phi = \sum_{\ell=1}^{L} J_\ell K
\)
denote the total number of random features across all levels. The parameter vector collects all coefficients
\(
\bm{a} = \{ a_{jk}^{(\ell)} \} \in \mathbb{R}^{N_\Phi}.
\)
Proceeding as in \cite{vanBeek2026RRQR}, the loss functional can be written as
\begin{equation}
\label{eq:multilevel_ls_all_constraints}
\min_{\bm{a} \in \mathbb{R}^{N_\Phi}}
\|\mat{N}\bm{a} - \mathbf{f}\|^2
+
\|\mat{B}\bm{a} - \mathbf{g}\|^2,
\end{equation}
where $\mat{N}$ and $\mat{B}$ denote the system matrices associated with the PDE residuals and the boundary conditions, respectively. The vectors $\mathbf{f}$ and $\mathbf{g}$ contain evaluations of the forcing term and boundary terms at interior and boundary collocation points, respectively.

As in the single-level formulation, the two residual contributions can be combined into a single global least-squares problem
\begin{equation}
\label{eq:multilevel_ls}
\min_{\bm{a} \in \mathbb{R}^{N_\Phi}}
\|\mat{M}\bm{a} - \bh\|^2,\quad\text{ where}\quad
\mat{M} =
\begin{bmatrix}
\mat{N} \\
\mat{B}
\end{bmatrix},
\qquad
\bh =
\begin{bmatrix}
\mathbf{f} \\
\mathbf{g}
\end{bmatrix}.
\end{equation}

\paragraph{Remark on linearity and extension to nonlinear PDEs.}
The formulation above relies critically on the assumption that the differential operator $\mathcal{N}$ and boundary operators $\mathcal{B}_k$ are linear, which enables the reduction of the training problem to a global least-squares system. This assumption, however, restricts the applicability of the current method to linear PDEs. For nonlinear problems, the direct least-squares formulation is no longer available. One possible extension would be to employ iterative linearisation strategies, such as Newton or Gauss--Newton methods, where each iteration involves solving a linearised least-squares problem. We note that such extensions would partially reintroduce iterative procedures but depending on the problem - like in the case of Newton-Krylov methods -  resulting linear systems are likely to be related and similar ``tricks" can be reused to regain efficiency. A systematic investigation of nonlinear extensions is left for future work.

\paragraph{Indexing conventions}
To describe the structure of the system matrices, we introduce the following indexing conventions.

\begin{itemize}
\item Let
\(
q = \sum_{m=1}^{\ell-1}J_m K + (j-1)K + k
\)
map the triple $(\ell,j,k)$ to a global column index.
\item Let
\(
p = \sum_{m=1}^{b-1} N_B^{(m)} + i
\)
map the pair $(i,b)$ to a global boundary collocation index.
\end{itemize}
We also denote
\(
N_B = \sum_{b=1}^{B} N_B^{(b)}
\)
as the total number of boundary collocation points.

\paragraph{Definition of vectors}
The vectors appearing in \eqref{eq:multilevel_ls} are defined as
\[
\mathbf{f} =
\left\{
\frac{1}{\sqrt{N_I}} f(\bx_n)
\right\}_{n=1}^{N_I},\qquad 
\mathbf{g} =
\{ g_p \}_{p=1}^{N_B},
\qquad
g_p =
\sqrt{\frac{\lambda_b}{N_B^{(b)}}}
g_b(\bx_i^{(b)}).
\]

\paragraph{Definition of the system matrices}
The matrices $\mat{N} \in \mathbb{R}^{N_I \times N_\Phi}$ and $\mat{B} \in \mathbb{R}^{N_B \times N_\Phi}$ are defined componentwise by
\[
N_{n,q}
=
\frac{1}{\sqrt{N_I} L}
\mathcal{N}
\left[
(\omega_j^{(\ell)} \phi_{jk}^{(\ell)})(\bx_n)
\right],\qquad
B_{p,q}
=
\sqrt{\frac{\lambda_b}{N_B^{(b)}}}\frac{1}{L}
\mathcal{B}_b
\left[
(\omega_j^{(\ell)} \phi_{jk}^{(\ell)})(\bx_i^{(b)})
\right].
\]

Each column therefore corresponds to the action of the PDE or boundary operator applied to a localised random feature
\(
\omega_j^{(\ell)}(\bx)\phi_{jk}^{(\ell)}(\bx).
\)
\paragraph{Multilevel domain-decomposition structure}
Because the window functions $\omega_j^{(\ell)}$ have compact support, the resulting matrix
\[
\mat{M} \in \mathbb{R}^{(N_I + N_B) \times N_\Phi}
\]
is sparse. In particular, an entry $M_{r,q}$ is nonzero only when
\(
\omega_j^{(\ell)}(\bx_r) \neq 0.
\)
We therefore define the index sets
\[
I_j^{(\ell)} =
\{ r \mid \omega_j^{(\ell)}(\bx_r) \neq 0 \},\quad
K_j^{(\ell)} =
\{ \sum_{m=1}^{\ell-1}J_m K + (j-1)K + 1,\dots,\sum_{m=1}^{\ell-1}J_m K + jK \},
\]
where $I_j^{(\ell)}$ contains the rows corresponding to collocation points influenced by subdomain $j$ at level $\ell$, and $K_j^{(\ell)}$ contains the associated feature indices. Let
\[
\mat{M}_j^{(\ell)}
=
\mat{M}|_{I_j^{(\ell)} \times K_j^{(\ell)}}
\]
denote the corresponding submatrix. Introducing restriction matrices
\[
\mat{V}_j^{(\ell)} \in
\mathbb{R}^{|I_j^{(\ell)}| \times (N_I+N_B)},
\qquad
\mat{W}_j^{(\ell)} \in
\mathbb{R}^{K \times N_\Phi},
\]
we can write
\[
\mat{M}_j^{(\ell)}
=
\mat{V}_j^{(\ell)} \mat{M} \mat{W}_j^{(\ell)T}.
\]
This yields the multilevel decomposition
\begin{equation}
\label{eq:multilevel_domain_decomp}
\mat{M}
=
\sum_{\ell=1}^{L}
\sum_{j=1}^{J_\ell}
\mat{V}_j^{(\ell)T}
\mat{M}_j^{(\ell)}
\mat{W}_j^{(\ell)} .
\end{equation}
Equation \eqref{eq:multilevel_domain_decomp} reveals the hierarchical block-sparse structure of the system matrix induced by the multilevel domain decomposition. An example plot of this sparse structure is shown in \Cref{fig:ch3-harm_osc_p0_preds}. In practice, this structure enables efficient assembly of $\mat{M}$ using sparse data structures and efficient solution of \eqref{eq:multilevel_ls} using iterative least-squares solvers such as LSQR. Once the coefficient vector $\bm{a}$ has been computed, the approximate solution $\hat{u}$ can be evaluated at arbitrary points by assembling the corresponding localised random feature expansions across all levels.

\section{Numerical assessment of the ELM-FBPINN method}
\label{sec:ch3_problems}
Our objective is to assess whether replacing gradient-based subdomain training by a random-feature least-squares solve can improve the efficiency and robustness of physics-informed solvers without \rev{increasing the error in the solution.} 
To isolate the effect of the training strategy from the effect of domain localisation, we compare three models:
(i) a standard single-domain PINN trained by gradient descent, 
(ii) an FBPINN that uses a domain decomposition and partition-of-unity construction and trainable subdomain networks, and 
(iii) the proposed ELM-FBPINN, which retains the FBPINN 
domain decomposition and 
partition-of-unity structure while replacing each subdomain network by a random-feature model and solving for the output coefficients via least squares (\Cref{tab:model-comparison}). We also assess the performance impact of including multiple levels in the ELM-FBPINN.
\begin{table}[t]
\centering
\small
\setlength{\tabcolsep}{6pt}
\renewcommand{\arraystretch}{1.25}
\begin{tabular}{lccc}
\toprule
 & \textbf{PINN} & \textbf{FBPINN} & \textbf{ELM-FBPINN} \\
\midrule
Domain decomposition / POU  & \ding{55} & \ding{51} & \ding{51} \\
Local (subdomain) models    & \ding{55} & \ding{51} & \ding{51} \\
Hidden features trainable?  & \ding{51} & \ding{51} & \ding{55} \\
Training strategy           & GD (Adam)  & GD (Adam)  & LS solver (e.g.\ LSQR) \\
Unknowns solved for         & all weights & all weights & output weights only \\
Linear solve required?      & \ding{55} & \ding{55} & \ding{51} \\
Expected effect             & spectral bias & mitigates bias & mitigates bias + faster solve \\
\bottomrule
\end{tabular}
\caption{Summary of the three compared approaches. }
\label{tab:model-comparison}
\end{table}
\rev{Three} representative test cases are considered:
\begin{itemize}
    \item \textbf{1D damped harmonic oscillator.} 
    This problem yields oscillatory solutions with controlled frequency content and exponential decay. It provides a controlled setting to probe spectral bias, the benefit of domain localisation, and the scaling of solver performance as the oscillation frequency increases.
    \item \textbf{2D multi-scale Laplacian problem.}
    This elliptic PDE admits an analytic solution composed of multiple sine modes, allowing systematic control of multiscale complexity in a higher-dimensional setting. 
    \rev{\item \textbf{2D inhomogeneous Helmholtz equation problem.}
    This problem is defined on an L-shaped domain with spatially varying Dirichlet boundary conditions and a non-trivial manufactured solution with a frequency parameter which allows systematic control of oscillatory complexity and optimisation difficulty.}
\end{itemize}
The numerical study is designed to isolate the influence of a small number of experimental parameters corresponding to the research questions introduced in Section~\ref{sec:introduction}. In particular, we vary the following levers:
\begin{itemize}
\item \textbf{RQ1: Training strategy}-gradient descent (GD) versus linear least-squares (LSQ).
\item \textbf{RQ2: Domain localisation}-number of subdomains $J_l$.
\item \textbf{RQ3: Basis enrichment}-number of basis functions per subdomain $K$.
\item \textbf{RQ4: Problem difficulty}-frequency or multiscale content scaled with $J_l$, and assessing the impact of including multiple levels, $L$.
\end{itemize}
{Although it is also possible to alter the depth ($h$) of the networks used, preliminary experiments indicated that this did not offer much benefit, and that shallow networks were optimal or near-optimal.}
The comparisons in \Cref{tab:model-comparison} and the research questions define our experimental framework. By separating localisation (shared by FBPINN and ELM-FBPINN) from training strategy (gradient descent versus least-squares), and by systematically varying capacity, and problem difficulty, the numerical study is designed to disentangle approximation effects from optimisation effects. 

\rev{
\paragraph{Remark on experimental scope.}
The experimental study focuses on comparisons with PINN and FBPINN in order to isolate the effect of the proposed ELM-based formulation within a controlled setting. While this enables a clear assessment of the methodological contribution, it does not aim to provide an exhaustive benchmark against all existing approaches. Extending the comparison to other scientific machine learning methods and more complex PDE settings is an important direction for future work.
}

\subsection{Common implementation details}
\label{sec:ch3_common_impl}
We now describe the common implementation choices used across all experiments which also represent a practical compromise as there are many possible choices of hyperparameters. These details ensure that comparisons are fair and that observed differences can be attributed to structural properties of the methods rather than to incidental training or sampling choices. Several aspects of the training and evaluation setup are shared across all experiments in one and two dimensions. Unless stated otherwise, the settings introduced below are fixed throughout all numerical experiments.

\paragraph{Domain decomposition} Both the FBPINN and ELM-FBPINN employ the same domain decomposition strategy. We assume that the global computational domain has a tensor-product structure,
\[
\Omega = \bigotimes_{i=1}^{d} [x^i_{\min}, x^i_{\max}],
\]
and at each level, $\ell$, we decompose $\Omega$ into overlapping equally sized intervals in one dimension and overlapping rectangular subdomains arranged on a uniform grid in two dimensions. Example domain decompositions in one and two dimensions are shown in \Cref{fig:ch3-harm_osc_p0_preds,fig:ch3-laplace-p0-preds}. \rev{When studying a two dimensional L-shaped problem domain in \Cref{sec:ch3_num_res_helmholtz}, we discard rectangular subdomains outside of the L-shaped domain, as shown in \cref{fig:ch3-helmholtz-p0-preds}.}

Localisation is enforced through a partition-of-unity (POU) construction at each level. The window functions are defined as
\begin{equation}
\label{eq:window_func}
\omega_{j}^{(\ell)}(\bx) =
\frac{\hat{\omega}_{j}^{(\ell)}(\bx)}{\sum_{j=1}^{J_{\ell}} \hat{\omega}_{j}^{(\ell)}(\bx)},
\qquad
\hat{\omega}_{j}^{(\ell)}(\bx) =
\prod_{i=1}^{d}
\Biggl[
1 + \cos\Bigl( \pi \frac{x_i - \mu^{(\ell)}_{j_i}}{\sigma^{(\ell)}_{j_i}} \Bigr)
\Biggr]^2,
\end{equation}
where 
$\mu^{(\ell)}_{j_i}$
and
$\sigma^{(\ell)}_{j_i}$ denote the centre and half-width of the $j$th subdomain in dimension $i$ at level $\ell$. 
For a one-dimensional domain, the subdomain bounds are defined as 
\begin{equation}
x^{(\ell)}_{j, \text{min}} = \mu^{(\ell)}_j - \frac{\delta}{2} (\mu^{(\ell)}_j - \mu^{(\ell)}_{j-1}),
\quad
x^{(\ell)}_{j, \text{max}} = \mu^{(\ell)}_j + \frac{\delta}{2} (\mu^{(\ell)}_{j+1} - \mu^{(\ell)}_j),
\end{equation}
with $\mu^{(\ell)}_j$ equally spaced over the domain and $\sigma^{(\ell)}_{j_i} = (x^{(\ell)}_{j, \text{max}} - x^{(\ell)}_{j, \text{min}})/2$, where $\delta$ is a width ratio that controls the degree of overlap and therefore information exchange between subdomains. An analogous definition is used for two dimensional domains with the same $\delta$ along each dimension. Note that a width ratio of less than one means that the subdomains no longer overlap. 
This cosine-based construction yields smooth window functions with compact support, ensuring locality while preserving a smooth global approximation.

\paragraph{Neural network architecture and weight initialisation} All models employ fully connected neural networks with $\tanh$ activation functions. In the PINN and FBPINN models, all network parameters are trainable. In contrast, the ELM-FBPINN fixes the hidden-layer parameters after random initialisation and trains only the linear output coefficients. All network weights are initialised using uniform LeCun initialisation \citep{lecun_efficient_1998}, i.e. randomly sampled according to $W_j^{(\ell,m)},b_j^{(\ell,m)} \sim \mathcal{U}(-a,a)$, where $a = R\sqrt{\frac{1}{\mathrm{dim}(\mathbf{z}^{\left(\ell,m-1\right)})}}$. Here $R$ is a positive parameter which we use to test the sensitivity of our models to the magnitude of the initialized weights.

\paragraph{Optimisation and solvers} 
The PINN and FBPINN models are trained using the Adam optimiser with a learning rate of $10^{-3}$ for $15,000$ iterations. The ELM-FBPINN replaces gradient-based optimisation with a sparse linear least-squares solve, which is performed using the iterative LSQR solver from SciPy \citep{virtanen_scipy_2020}, 
with a maximum of $15,000$ iterations. Unless stated otherwise, the solver is allowed to run for the full iteration budget in order to \rev{facilitate %direct 
comparison of convergence behaviour. We note, however, that one iteration of LSQR may not be equivalent to one iteration of Adam}. In practice, many problems studied here converge much faster than the maximum number iterations and could be stopped early with an appropriate stopping condition. Model performance is assessed using the \rev{test error, measured in the $e_{L^1}$ norm,}  
given by
\begin{equation}
\label{l1_error}
e_{L^1} =
\frac{1}{M \varsigma}\sum_{i=1}^{M }
\lvert \hat{u}(\mathbf{x}_i) - u(\mathbf{x}_i) \rvert,
\end{equation}
where $M$ is the number of test points, and $\varsigma$ is the standard deviation of the set of true solutions $\{u(\mathbf{x}_i)\}_i^M$. Training and test points are sampled on a regular grid, with sufficient resolution per subdomain to accurately capture local approximation quality. For fairness exactly the same collocation points and random seeds are used for each model compared (PINN/FBPINN/ELM-FBPINN). For robustness each model is retrained using five independent random seeds, and test errors are reported as the median value, with error bars indicating the range between the best and worst observed values. Training times are reported as the mean with error bars indicating one standard deviation. For the ELM-FBPINN, training time includes every stage of the computation; construction of the local bases, assembly of the global system, and iterative solution time.

\paragraph{Software and hardware implementation} All models are implemented in Python using JAX \cite{bradbury_jax_2018}, building upon the highly optimised open-source FBPINN library \cite{Moseley2023}. All models are trained on a single AMD EPYC 7742 CPU core with 50 GB RAM, using the SciPy library \cite{2020SciPy-NMeth} for defining sparse linear solvers.

%% FIGURE CAPTION MACROS

\newcommand{\captionfigbase}[2]{%
\caption{
\rev{Test error, measured in the $e_{L^1}$ norm,} versus iteration count for the baseline #1 problem, where $J$ lists the number of subdomains #2 each level (here, in all cases $L=1$). \rev{The median test error is denoted by a solid line, while the shaded region indicates the range between the best and worst seeds.}
}%
}

\newcommand{\captiontabbase}[2]{%
\caption{
Baseline results for the #1 problem, where $J$ lists the number of subdomains #2 each level (here, in all cases $L=1$). \rev{The values for the test error, measured in the $e_{L^1}$ norm, state the median value together with the range between the best and worst observed values. For runtime measurements, we report the mean together with one standard deviation.}
}%
}

\newcommand{\captionfigK}[2]{%
\caption{
\rev{Test error, measured in the $e_{L^1}$ norm,} vs number of basis functions, $K$, for the #1 problem, where $J$ lists the number of subdomains #2 each level (here, in all cases $L=1$). \rev{In the left-hand plot, the median test error is denoted by a solid line, while the shaded region indicates the range between the best and worst seeds. In the right-hand plot, the median final test error is denoted by a circle, and the range between the best and worst seeds by a vertical bar.}
}%
}

\newcommand{\captionfigJ}[2]{%
\caption{
\rev{Test error, measured in the $e_{L^1}$ norm,} vs number of subdomains, $J$, for FBPINN and ELM-FBPINN, for the #1 problem, where $J$ lists the number of subdomains #2 each level (here, in all cases $L=1$). \rev{In the left-hand plot, the median test error is denoted by a solid line, while the shaded region indicates the range between the best and worst seeds. In the right-hand plot, the median final test error is denoted by a circle, and the range between the best and worst seeds by a vertical bar.}
}%
}

\newcommand{\captionfigscaling}[3]{%
\caption{
\rev{Test error, measured in the $e_{L^1}$ norm,} as the #1 increases (with $J$ scaled accordingly) for the #2 problem, where $J$ lists the number of subdomains #3 each level. \rev{In the left-hand plot, the median test error is denoted by a solid line, while the shaded region indicates the range between the best and worst seeds. In the right-hand plot, the median final test error is denoted by a circle, and the range between the best and worst seeds by a vertical bar.}
}%
}

\newcommand{\captiontabK}[2]{%
\caption{
Test error and runtime when varying number of basis functions, $K$, for the #1 problem, where $J$ lists the number of subdomains #2 each level (here, in all cases $L=1$). \rev{The values for the test error, measured in the $e_{L^1}$ norm, state the median value together with the range between the best and worst observed values. For runtime measurements, we report the mean together with one standard deviation.}
}%
}

\newcommand{\captiontabJ}[2]{%
\caption{
Test error and runtime when varying number of subdomains, $J$, for the #1 problem, where $J$ lists the number of subdomains #2 each level (here, in all cases $L=1$). \rev{The values for the test error, measured in the $e_{L^1}$ norm, state the median value together with the range between the best and worst observed values. For runtime measurements, we report the mean together with one standard deviation.}
}%
}

\newcommand{\captiontabscaling}[3]{%
\caption{
Test error and runtime when varying #1 for the #2 problem, where $J$ lists the number of subdomains #3 each level. \rev{The values for the test error, measured in the $e_{L^1}$ norm, state the median value together with the range between the best and worst observed values. For runtime measurements, we report the mean together with one standard deviation.}
}%
}

\subsection{1D problem -- Damped harmonic oscillator}
\label{sec:ch3_problem_harmosc}

The damped harmonic oscillator provides a controlled setting in which the solution exhibits oscillations with tunable frequency and exponential decay. This setting probes robustness to increasing frequency, the benefit of domain decomposition, and the effect of replacing nonlinear optimisation by a linear least-squares solve.

\paragraph{Problem definition}
We consider the displacement $u(t)$ on $t\in(0,1]$ governed by
\begin{equation}
\label{eq:harm_oscillator}
\left\{
\begin{array}{l}
m u_{tt} + \mu u_t + k u(t) = 0, \quad t \in (0,1], \\
u(0) = u_0, \quad u_{t}(0) = v_0,
\end{array}
\right.
\end{equation}
where $m$ is the mass, $\mu$ the damping coefficient, and $k$ the spring constant.

\paragraph{Under-damped regime and tunable frequency}
Defining
\[
{\dho} = \frac{\mu}{2m},
\qquad
\omega_0 = \sqrt{\frac{k}{m}},
\]
the under-damped regime corresponds to ${\dho} < \omega_0$, with oscillation frequency
\[
\omega = \sqrt{\omega_0^2 - {\dho}^2}.
\]
We fix $u_0=1$, $v_0=0$, $m=1$ and $\mu=4$ (hence $\dho=2$) and vary $\omega_0$ to control difficulty: increasing $\omega_0$ introduces higher-frequency components while preserving the decay rate. The exact solution is
\begin{equation}
\label{eq:harm_u_exact}
u(t) = e^{-{\dho} t}\, 2A \cos(\phi + \omega t),
\end{equation}
with
\[
A = \frac{u_0}{2\cos(\phi)},
\qquad
\phi = \arctan\!\Big(-\frac{(v_0/u_0)+{\dho}}{\omega}\Big).
\]

\paragraph{Physics-informed objective}
All models minimise the same residual-based loss,
\begin{equation}
\mathcal{L}
=
\mathcal{L}_{\mathrm{phys}}
+
\lambda_{u_0}\mathcal{L}_{\mathrm{ic}}^{(u_0)}
+
\lambda_{v_0}\mathcal{L}_{\mathrm{ic}}^{(v_0)},
\end{equation}
with
\begin{align}
\mathcal{L}_{\mathrm{phys}}
&=
\frac{1}{N_I}
\sum_{i=1}^{N_I}
\big(
m\hat{u}_{tt} + \mu\hat{u}_t + k\hat{u}
\big)^2(t_i), \quad t_i \in (0,1], \\
\mathcal{L}_{\mathrm{ic}}^{(u_0)}
&=
\big(\hat{u}(0)-u_0\big)^2, \\
\mathcal{L}_{\mathrm{ic}}^{(v_0)}
&=
\big(\hat{u}_t(0)-v_0\big)^2.
\end{align}
\rev{We set $\lambda_{u_0}=k^2=\omega_0^4$ and $\lambda_{v_0}=k=\omega_0^2$, which normalises all loss terms to the same physical unit and we find provides a reasonable balance between them as the solution frequency increases.}

\paragraph{Summary of observed behaviour}
On highly oscillatory instances, a standard PINN fails to converge within the training budget, whereas both FBPINN and ELM-FBPINN converge reliably. Across ablations, ELM-FBPINN attains comparable or lower errors than FBPINN while converging faster \rev{in terms of number of iteration steps and wall-clock time}, provided each subdomain has sufficient local capacity. Failure modes are dominated by insufficient localisation (small $J_l$) or feature degeneracy (small $K$), rather than optimisation speed alone.

\paragraph{Ablation roadmap}
Starting from a baseline configuration, we vary (i) local capacity ($K$), (ii) localisation resolution ($J_l$), and (iii) frequency scaling via coupled $\omega_0$ and $J_l$, comparing one-level with multilevel models. When relevant, we report the condition number and numerical rank of the ELM-FBPINN system matrix to interpret stability trends.

% ============================================================
\subsubsection{Baseline}
\label{sec:ch3_harm_baseline}

We consider the under-damped oscillator with $\omega_0 = 80$,
yielding a highly oscillatory solution with slowly decaying amplitude.
The PINN uses $h=2$ hidden layers and $K=64$ basis functions per layer.
Both the FBPINN and ELM-FBPINN use a single hidden layer ($h=1$) with $K=12$ basis functions per subdomain, and a one-level decomposition into $J_1=20$ subdomains with width ratio $\delta=2.9$. \rev{The PINN and FBPINN use a fully connected network (FCN) architecture, while the ELM-FBPINN uses an extreme learning machine (ELM) architecture.} All models use a weight scaling of $R=1$ and are trained on 260 regularly spaced collocation points.

\begin{figure}[h]
    \centering
    \includegraphics[width=0.9\linewidth]{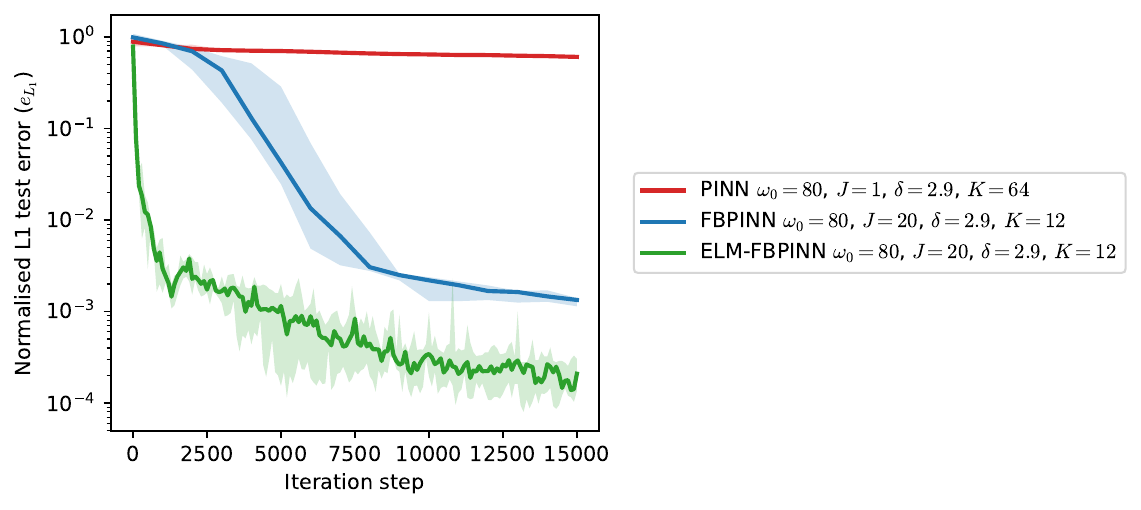}
\captionfigbase{1D damped harmonic oscillator}{for}
    \label{fig:ch3-harm-baseline-l1-iterations}
\end{figure}

\Cref{fig:ch3-harm-baseline-l1-iterations} shows the \rev{test error, measured in the $e_{L^1}$ norm,} with $M=1000$ regularly spaced test points versus iteration count.
The PINN does not converge to an acceptable solution within the training budget, consistent with the well-known difficulty of gradient-based PINNs in resolving high-frequency components.
Both decomposed methods converge, and ELM-FBPINN reaches the lowest final error while converging faster \rev{in terms of number of iteration steps} than FBPINN. While PINN fails to capture the real solution leading to large errors, the two decomposed methods yield similar predictions. \Cref{tab:ch3-harm-p0-results} summarises final test errors and wall-clock times. In this baseline, ELM-FBPINN attains the best median \rev{test error, measured in the $e_{L^1}$ norm,} with a faster training time to the PINN and FBPINN.

\begin{figure}[h]
    \centering
    \includegraphics[width=0.9\linewidth]{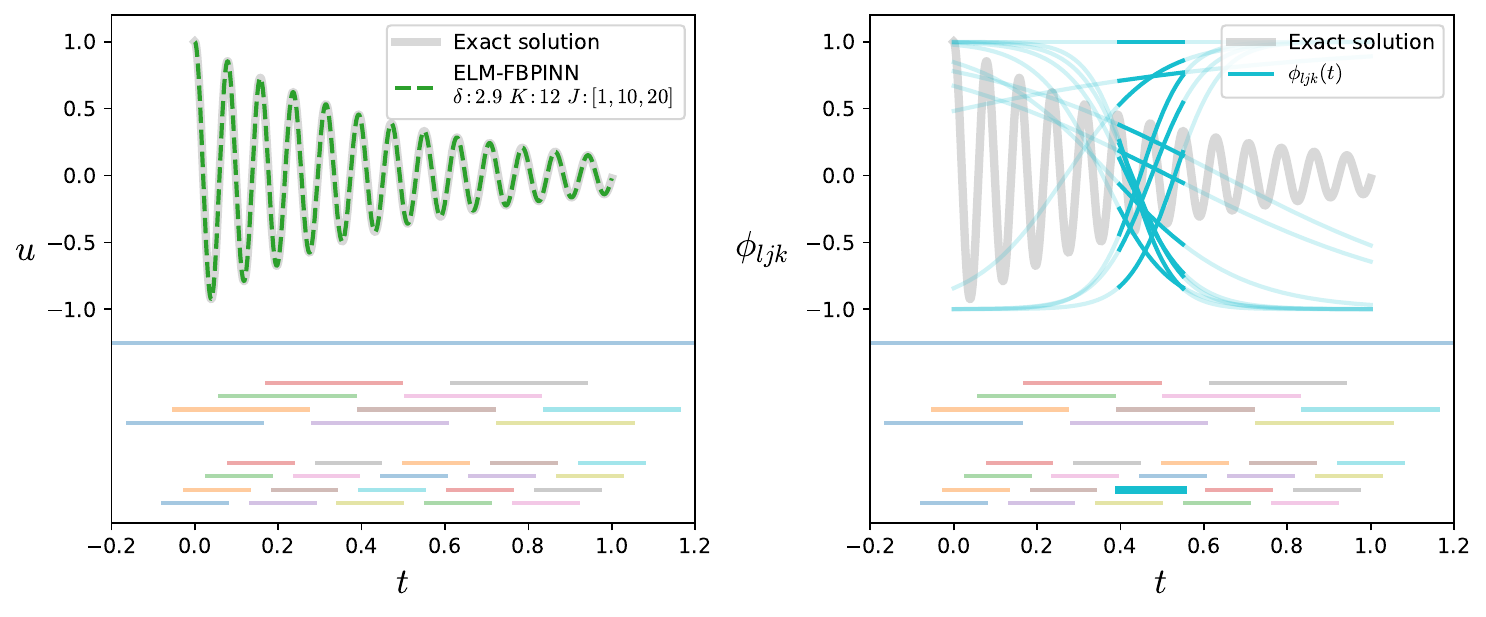}
    \vspace{0.5cm}
    \includegraphics[width=0.45\linewidth]{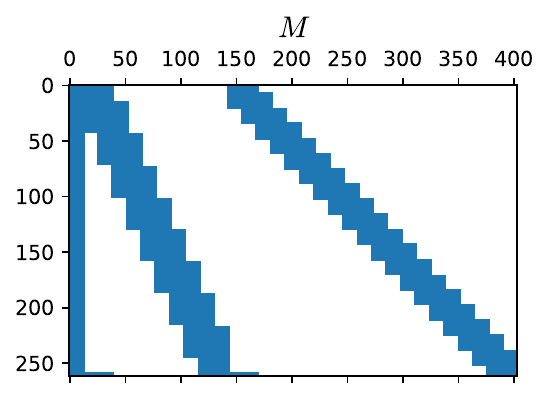}
    \caption{Model predictions versus the exact solution for the 1D damped harmonic oscillator with a multilevel ELM-FBPINN, where $J$ lists the number of subdomains for each level; basis functions generated by randomisation; sparse multilevel structure of \rev{the system matrix $\mat{M}$ in \cref{eq:multilevel_domain_decomp}.}}
    \label{fig:ch3-harm_osc_p0_preds}
\end{figure}

\begin{table}[h]
\setlength{\tabcolsep}{2pt}
\centering
\resizebox{\textwidth}{!}{
\begin{tabular}{cccccccc}
\toprule
Model & $(\omega_{0},J)$ & $(h,K)$ & $\delta$ & $\kappa(M)$ & Optimiser & $e_{L_1}$ & Time (s) \\
\midrule
PINN & (80, 1) & (2, 64) & 2.9 & N/A & Adam & 6.1$\pm$0.4e-01 & 2.9$\pm$0.0e+01 \\
FBPINN & (80, 20) & (1, 12) & 2.9 & N/A & Adam & 1.3$\pm$0.3e-03 & 1.9$\pm$0.0e+01 \\
ELM-FBPINN & (80, 20) & (1, 12) & 2.9 & 8.5e+18 & LSQR & 1.1$\pm$0.5e-04 & 5.4$\pm$0.0 \\
\bottomrule
\end{tabular}

}
\captiontabbase{1D damped harmonic oscillator}{for}
\label{tab:ch3-harm-p0-results}
\end{table}

We also show the model predictions in \Cref{fig:ch3-harm_osc_p0_preds} for the same ELM-FBPINN, but with $L=3$ and  $J = [1, 10, 20]$ subdomains, where $J$ lists the number of subdomains for each level (i.e., $J_1=1, J_2=10$, and $J_3=20$). The multilevel solution is visually indistinguishable from the exact solution. The figure also illustrates the randomly generated basis functions and the sparse multilevel structure of the system matrix \rev{$\mat{M}$}, which is a key feature of the proposed approach.

% ============================================================
\subsubsection{Varying number of basis functions}
\label{sec:ch3_harm_ablation_C}

\begin{figure}[h]
    \centering
    \includegraphics[width=\textwidth]{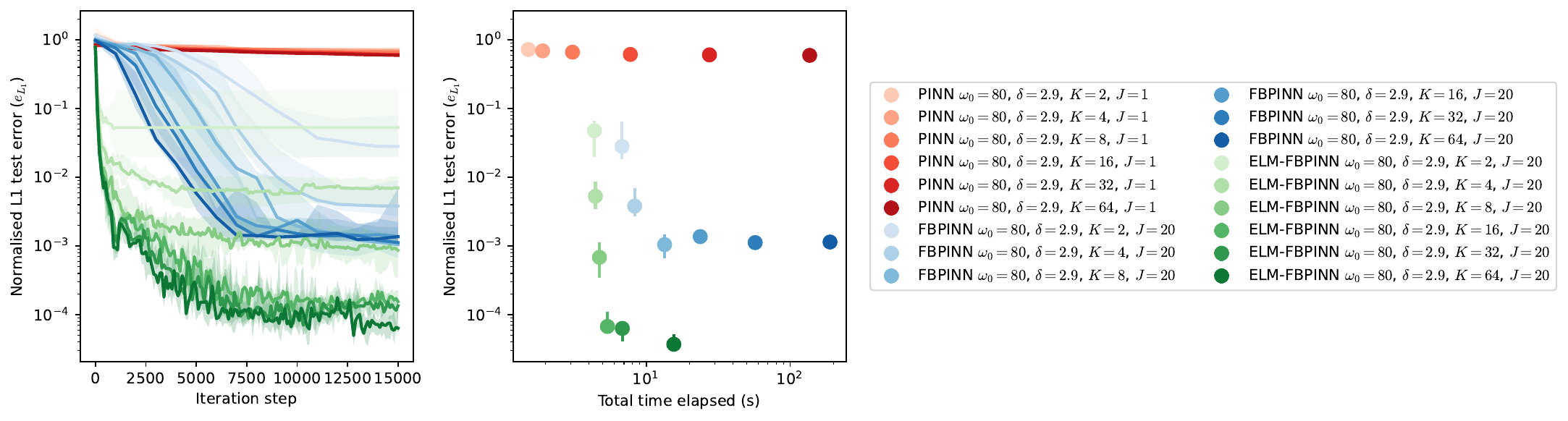}
\captionfigK{1D damped harmonic oscillator}{for}
    \label{fig:ch3-harm-c_ablation-final_l1}
\end{figure}

We vary the number of basis functions $K$ of the baseline $L=1$ model to probe how each method exploits increased local capacity.
This ablation distinguishes (i) expressivity limitations at small $K$ from (ii) redundancy/feature correlation at large $K$. Convergence curves and final errors are shown in \Cref{fig:ch3-harm-c_ablation-final_l1} and \Cref{tab:ch3-harm-c-ablation}. \rev{Note, for the baseline PINN which uses $h=2$, varying $K$ also varies the number of basis functions per hidden layer as we use $K$ basis functions per hidden layer.}

\paragraph{Key observations}
(i) All PINN configurations fail to converge, even as $K$ increases, reinforcing that additional capacity alone does not resolve the optimisation difficulty on this high-frequency problem.
(ii) Both decomposed methods converge for moderate $K$, but ELM-FBPINN is the only method that shows a strong and systematic improvement with increasing $K$ before saturating.
(iii) At very small $K$, FBPINN can outperform ELM-FBPINN: in this low-capacity regime, the random-feature trial space is too limited to generalise accurately, and the linear solve cannot compensate for missing expressivity.

For the ELM-FBPINN, the error decreases rapidly as the number of basis functions $K$ increases from very small values, and then plateaus for larger $K$, indicating that once the dominant oscillatory structure is resolved, additional basis functions provide diminishing returns.
This occurs despite a rapid increase in the condition number $\kappa(M)$, suggesting growing linear dependence among the random basis functions, although the least-squares solver remains stable in practice.
In contrast, the PINN shows only marginal improvement as $K$ increases while the computational cost grows significantly, whereas the FBPINN reaches good \rev{test error} for moderate $K$ but also plateaus.
Overall, the ELM-FBPINN achieves the lowest errors across all tested values of $K$ while maintaining faster runtimes, since training reduces to the solution of a linear least-squares problem rather than iterative gradient-based optimisation.

\begin{table}[h]
\setlength{\tabcolsep}{2pt}
\centering
\resizebox{\textwidth}{!}{
\begin{tabular}{cccccccc}
\toprule
Model & $(\omega_{0},J)$ & $(h,K)$ & $\delta$ & $\kappa(M)$ & Optimiser & $e_{L_1}$ & Time (s) \\
\midrule
PINN & (80, 1) & (2, 2) & 2.9 & N/A & Adam & 7.2$\pm$0.1e-01 & 1.5$\pm$0.0 \\
PINN & (80, 1) & (2, 4) & 2.9 & N/A & Adam & 6.9$\pm$0.5e-01 & 1.9$\pm$0.0 \\
PINN & (80, 1) & (2, 8) & 2.9 & N/A & Adam & 6.6$\pm$0.2e-01 & 3.1$\pm$0.0 \\
PINN & (80, 1) & (2, 16) & 2.9 & N/A & Adam & 6.2$\pm$0.5e-01 & 7.8$\pm$0.0 \\
PINN & (80, 1) & (2, 32) & 2.9 & N/A & Adam & 6.1$\pm$0.4e-01 & 2.7$\pm$0.0e+01 \\
PINN & (80, 1) & (2, 64) & 2.9 & N/A & Adam & 6.0$\pm$0.4e-01 & 1.4$\pm$0.0e+02 \\
FBPINN & (80, 20) & (1, 2) & 2.9 & N/A & Adam & 2.8$\pm$4.7e-02 & 6.8$\pm$0.1 \\
FBPINN & (80, 20) & (1, 4) & 2.9 & N/A & Adam & 3.8$\pm$4.3e-03 & 8.3$\pm$0.2 \\
FBPINN & (80, 20) & (1, 8) & 2.9 & N/A & Adam & 1.0$\pm$0.8e-03 & 1.3$\pm$0.0e+01 \\
FBPINN & (80, 20) & (1, 16) & 2.9 & N/A & Adam & 1.4$\pm$0.6e-03 & 2.4$\pm$0.0e+01 \\
FBPINN & (80, 20) & (1, 32) & 2.9 & N/A & Adam & 1.1$\pm$0.3e-03 & 5.7$\pm$0.0e+01 \\
FBPINN & (80, 20) & (1, 64) & 2.9 & N/A & Adam & 1.1$\pm$0.3e-03 & 1.9$\pm$0.0e+02 \\
ELM-FBPINN & (80, 20) & (1, 2) & 2.9 & 1.1e+06 & LSQR & 4.8$\pm$4.6e-02 & 4.4$\pm$0.0 \\
ELM-FBPINN & (80, 20) & (1, 4) & 2.9 & 2.1e+09 & LSQR & 5.3$\pm$5.2e-03 & 4.4$\pm$0.0 \\
ELM-FBPINN & (80, 20) & (1, 8) & 2.9 & 1.6e+17 & LSQR & 6.8$\pm$7.9e-04 & 4.7$\pm$0.0 \\
ELM-FBPINN & (80, 20) & (1, 16) & 2.9 & 2.1e+19 & LSQR & 6.7$\pm$4.4e-05 & 5.4$\pm$0.0 \\
ELM-FBPINN & (80, 20) & (1, 32) & 2.9 & 2.0e+20 & LSQR & 6.3$\pm$2.8e-05 & 6.8$\pm$0.0 \\
ELM-FBPINN & (80, 20) & (1, 64) & 2.9 & 3.0e+20 & LSQR & 3.7$\pm$1.7e-05 & 1.6$\pm$0.0e+01 \\
\bottomrule
\end{tabular}

}
\captiontabK{1D damped harmonic oscillator}{for}
\label{tab:ch3-harm-c-ablation}
\end{table}

\paragraph{Interpretation via linear-system diagnostics}
As $K$ increases, $\kappa(\mathbf{M})$ grows rapidly and the numerical rank saturates, consistent with the introduction of increasingly correlated random features.
This provides a qualitative explanation for the observed \rev{plateau in the test error}: beyond moderate $K$, additional basis functions add limited independent information.

{Overall, except for very small $K$ ($K=2,4$), ELM-FBPINN outperforms both FBPINN and PINN in \rev{terms of the test error} while avoiding gradient-based optimisation. Despite the growth of $\kappa(\mathbf{M})$ with $K$, the least-squares problems are evidently solved well enough to obtain accurate ELM-FBPINN solutions. This suggests that the full least-squares system is not as numerically problematic as $\kappa(\mathbf{M})$ alone might indicate.}

% ============================================================
\subsubsection{Varying number of subdomains}
\label{sec:ch3_harm_ablation_J}

We vary the number of subdomains $J_1$  of the baseline $L=1$ model to test how localisation affects approximation quality and convergence.
Unlike increasing $K$, increasing $J_1$ primarily adds \emph{spatially distinct} basis functions, which should improve representation of highly oscillatory solutions once subdomains are sufficiently small.

\begin{figure}[t]
    \centering
    \includegraphics[width=\linewidth]{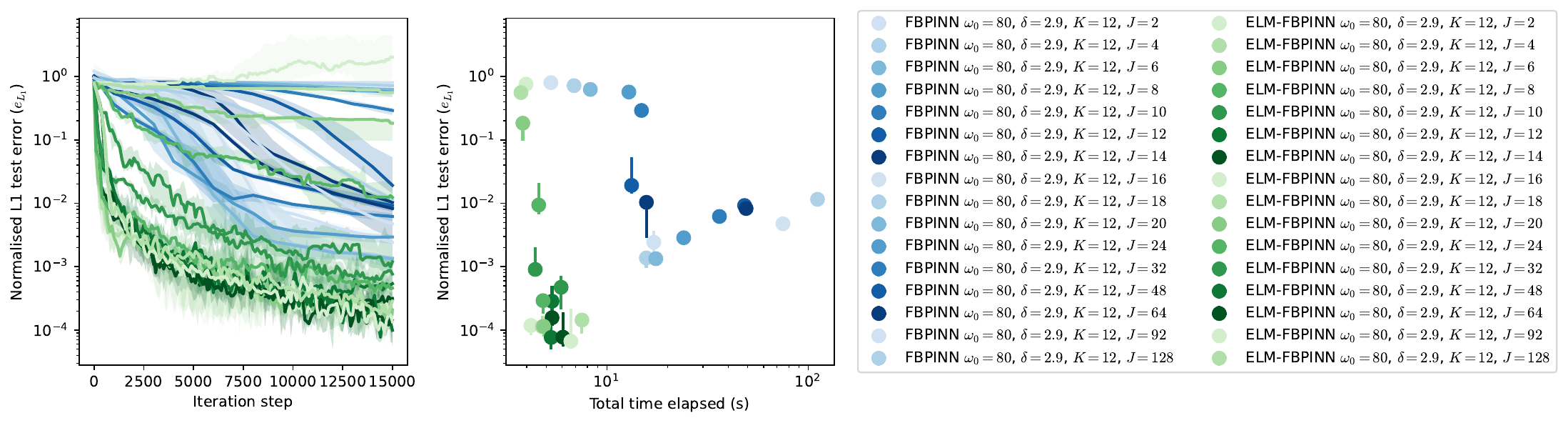}
\captionfigJ{1D damped harmonic oscillator}{for}
    \label{fig:ch3-harm-j-l1vsiter}
\end{figure}

\paragraph{Key observations}
\Cref{fig:ch3-harm-j-l1vsiter} shows that both methods improve with increasing $J_1$, but a minimum level of localisation is required: $J_1=2$ and $J_1=4$ fail to converge reliably for either approach.
For $J_1\gtrsim 8$, ELM-FBPINN converges faster \rev{in terms of number of iteration steps} than FBPINN and typically attains lower errors.
Final errors and times are summarised in \Cref{tab:ch3-harm-j-ablation}.
\rev{The test error} improves significantly as the number of subdomains $J_1$ increases from very small values, reflecting the benefit of domain localisation for resolving the oscillatory structure of the solution.
The error decreases rapidly up to moderate values of $J_1$ (approximately $J_1=12$--$20$), after which \rev{the test error} saturates and may slightly degrade for very large $J_1$, indicating diminishing returns once each local subproblem becomes sufficiently simple and the decomposition becomes excessively fine.
In contrast to increasing $K$, increasing $J_1$ does not dramatically worsen $\kappa(\mathbf{M})$, which remains in the range $10^{17}$--$10^{19}$, while the numerical rank grows approximately linearly with $J_1$, indicating that additional subdomains introduce largely independent degrees of freedom.
This supports the interpretation that increasing $J_1$ expands the trial space by adding spatially distinct basis functions rather than correlated features.
Finally, we observe that the computational cost grows significantly with $J_1$ for FBPINN due to the increased optimisation cost, whereas for ELM-FBPINN the runtime increases only mildly.

\begin{table}[h!]
\setlength{\tabcolsep}{2pt}
\centering
\resizebox{\textwidth}{!}{
\begin{tabular}{cccccccc}
\toprule
Model & $(\omega_{0},J)$ & $(h,K)$ & $\delta$ & $\kappa(M)$ & Optimiser & $e_{L_1}$ & Time (s) \\
\midrule
FBPINN & (80, 2) & (1, 12) & 2.9 & N/A & Adam & 8.0$\pm$0.0e-01 & 5.3$\pm$0.1 \\
FBPINN & (80, 4) & (1, 12) & 2.9 & N/A & Adam & 7.2$\pm$0.5e-01 & 6.9$\pm$0.0 \\
FBPINN & (80, 6) & (1, 12) & 2.9 & N/A & Adam & 6.3$\pm$0.8e-01 & 8.3$\pm$0.1 \\
FBPINN & (80, 8) & (1, 12) & 2.9 & N/A & Adam & 5.7$\pm$1.8e-01 & 1.3$\pm$0.0e+01 \\
FBPINN & (80, 10) & (1, 12) & 2.9 & N/A & Adam & 2.9$\pm$0.7e-01 & 1.5$\pm$0.0e+01 \\
FBPINN & (80, 12) & (1, 12) & 2.9 & N/A & Adam & 1.9$\pm$3.9e-02 & 1.3$\pm$0.0e+01 \\
FBPINN & (80, 14) & (1, 12) & 2.9 & N/A & Adam & 1.0$\pm$0.9e-02 & 1.6$\pm$0.0e+01 \\
FBPINN & (80, 16) & (1, 12) & 2.9 & N/A & Adam & 2.4$\pm$2.2e-03 & 1.7$\pm$0.0e+01 \\
FBPINN & (80, 18) & (1, 12) & 2.9 & N/A & Adam & 1.4$\pm$0.5e-03 & 1.6$\pm$0.0e+01 \\
FBPINN & (80, 20) & (1, 12) & 2.9 & N/A & Adam & 1.3$\pm$0.3e-03 & 1.7$\pm$0.0e+01 \\
FBPINN & (80, 24) & (1, 12) & 2.9 & N/A & Adam & 2.9$\pm$0.9e-03 & 2.4$\pm$0.0e+01 \\
FBPINN & (80, 32) & (1, 12) & 2.9 & N/A & Adam & 6.2$\pm$1.7e-03 & 3.6$\pm$0.0e+01 \\
FBPINN & (80, 48) & (1, 12) & 2.9 & N/A & Adam & 9.2$\pm$2.9e-03 & 4.8$\pm$0.0e+01 \\
FBPINN & (80, 64) & (1, 12) & 2.9 & N/A & Adam & 8.3$\pm$1.1e-03 & 4.9$\pm$0.1e+01 \\
FBPINN & (80, 92) & (1, 12) & 2.9 & N/A & Adam & 4.8$\pm$1.4e-03 & 7.5$\pm$0.1e+01 \\
FBPINN & (80, 128) & (1, 12) & 2.9 & N/A & Adam & 1.2$\pm$0.3e-02 & 1.1$\pm$0.0e+02 \\
ELM-FBPINN & (80, 2) & (1, 12) & 2.9 & 5.1e+17 & LSQR & 7.6$\pm$0.2e-01 & 4.0$\pm$0.0 \\
ELM-FBPINN & (80, 4) & (1, 12) & 2.9 & 1.7e+18 & LSQR & 5.6$\pm$1.3e-01 & 3.7$\pm$0.0 \\
ELM-FBPINN & (80, 6) & (1, 12) & 2.9 & 3.2e+18 & LSQR & 1.8$\pm$1.1e-01 & 3.8$\pm$0.0 \\
ELM-FBPINN & (80, 8) & (1, 12) & 2.9 & 4.9e+18 & LSQR & 9.4$\pm$14.4e-03 & 4.6$\pm$0.0 \\
ELM-FBPINN & (80, 10) & (1, 12) & 2.9 & 1.0e+18 & LSQR & 9.1$\pm$12.2e-04 & 4.4$\pm$0.0 \\
ELM-FBPINN & (80, 12) & (1, 12) & 2.9 & 2.2e+18 & LSQR & 2.9$\pm$3.3e-04 & 5.3$\pm$0.0 \\
ELM-FBPINN & (80, 14) & (1, 12) & 2.9 & 2.3e+18 & LSQR & 1.6$\pm$0.7e-04 & 5.3$\pm$0.0 \\
ELM-FBPINN & (80, 16) & (1, 12) & 2.9 & 4.5e+18 & LSQR & 1.2$\pm$0.5e-04 & 4.2$\pm$0.0 \\
ELM-FBPINN & (80, 18) & (1, 12) & 2.9 & 1.7e+19 & LSQR & 1.1$\pm$0.7e-04 & 4.8$\pm$0.0 \\
ELM-FBPINN & (80, 20) & (1, 12) & 2.9 & 8.5e+18 & LSQR & 1.1$\pm$0.5e-04 & 4.8$\pm$0.0 \\
ELM-FBPINN & (80, 24) & (1, 12) & 2.9 & 1.8e+19 & LSQR & 2.9$\pm$1.5e-04 & 4.8$\pm$0.1 \\
ELM-FBPINN & (80, 32) & (1, 12) & 2.9 & 1.5e+19 & LSQR & 4.7$\pm$5.0e-04 & 5.9$\pm$0.0 \\
ELM-FBPINN & (80, 48) & (1, 12) & 2.9 & 9.2e+18 & LSQR & 7.6$\pm$14.4e-05 & 5.3$\pm$0.0 \\
ELM-FBPINN & (80, 64) & (1, 12) & 2.9 & 1.4e+19 & LSQR & 7.7$\pm$13.5e-05 & 6.0$\pm$0.0 \\
ELM-FBPINN & (80, 92) & (1, 12) & 2.9 & 6.9e+19 & LSQR & 6.7$\pm$16.1e-05 & 6.6$\pm$0.1 \\
ELM-FBPINN & (80, 128) & (1, 12) & 2.9 & 8.6e+18 & LSQR & 1.4$\pm$0.9e-04 & 7.5$\pm$0.1 \\
\bottomrule
\end{tabular}

}
\captiontabJ{1D damped harmonic oscillator}{for}
\label{tab:ch3-harm-j-ablation}
\end{table}

% ============================================================
\subsubsection{Scaling solution frequency with subdomains}
\label{sec:ch3_harm_scaling_w0J}

 We test scalability to harder (higher-frequency) solutions by increasing $\omega_0$ while simultaneously increasing $J_1$ (and the number of collocation points proportionally) for the baseline $L=1$ model.
 This keeps the \emph{effective resolution per subdomain} approximately constant, so that changes in performance reflect the optimisation/solve strategy rather than a trivial loss of spatial resolution. In this section, we also assess whether including multiple levels can improve scaling performance by aiding global subdomain communication with the addition of coarser levels.

 \paragraph{Key observations}
 \Cref{fig:ch3-harm-multi-j-convergence,{tab:ch3-harm-j-multiscale}} show that both the FBPINN and ELM-FBPINN converge more slowly as $\omega_0$ increases, but ELM-FBPINN consistently achieves lower errors and shorter runtimes across all tested frequencies.
Moreover adding multiple levels significantly improves the \rev{test error} of the ELM-FBPINN, with only a small increase in training time. This suggests that adding corse levels aids global subdomain communication and facilitates faster convergence, whilst keeping the least squares system sparse and efficient to solve.

\begin{figure}[t]
    \centering
    \includegraphics[width=\linewidth]{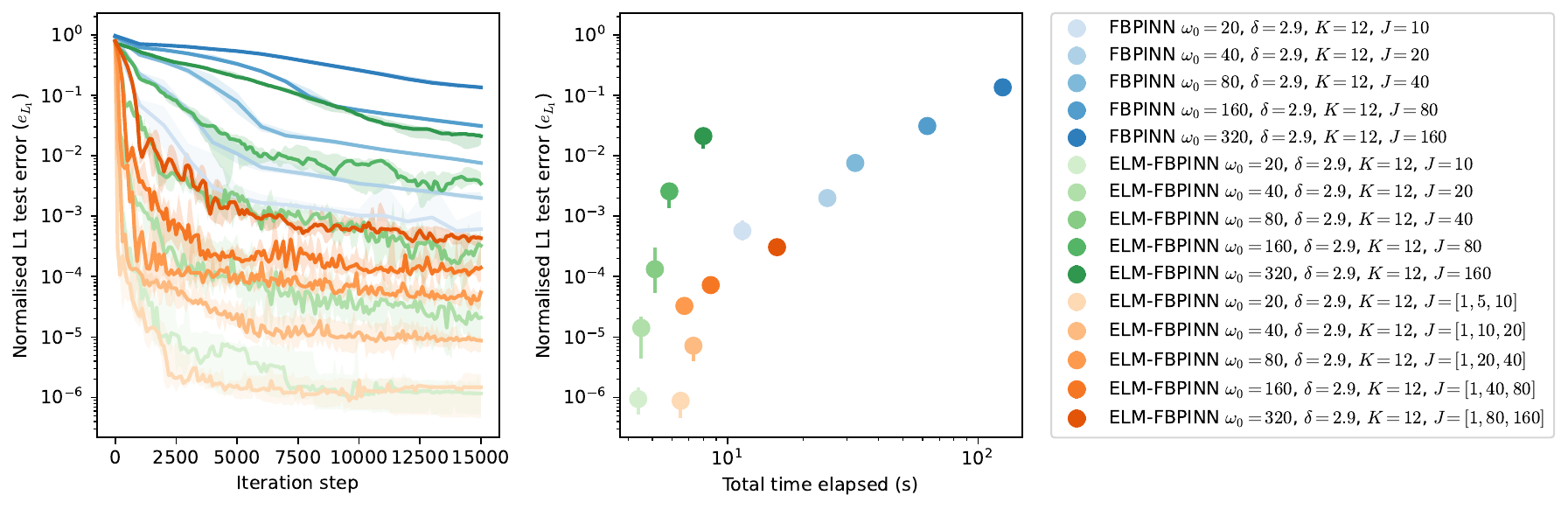}
\captionfigscaling{solution frequency $\omega_0$}{1D damped harmonic oscillator}{for}
    \label{fig:ch3-harm-multi-j-convergence}
\end{figure}

\begin{table}[h]
\setlength{\tabcolsep}{2pt}
\centering
\resizebox{\textwidth}{!}{
\begin{tabular}{cccccccc}
\toprule
Model & $(\omega_{0},J)$ & $(h,K)$ & $\delta$ & $\kappa(M)$ & Optimiser & $e_{L_1}$ & Time (s) \\
\midrule
FBPINN & (20, 10) & (1, 12) & 2.9 & N/A & Adam & 5.7$\pm$4.4e-04 & 1.1$\pm$0.0e+01 \\
FBPINN & (40, 20) & (1, 12) & 2.9 & N/A & Adam & 2.0$\pm$0.6e-03 & 2.5$\pm$0.1e+01 \\
FBPINN & (80, 40) & (1, 12) & 2.9 & N/A & Adam & 7.6$\pm$0.7e-03 & 3.2$\pm$0.0e+01 \\
FBPINN & (160, 80) & (1, 12) & 2.9 & N/A & Adam & 3.1$\pm$0.5e-02 & 6.3$\pm$0.2e+01 \\
FBPINN & (320, 160) & (1, 12) & 2.9 & N/A & Adam & 1.4$\pm$0.1e-01 & 1.3$\pm$0.0e+02 \\
ELM-FBPINN & (20, 10) & (1, 12) & 2.9 & 1.2e+18 & LSQR & 9.5$\pm$9.6e-07 & 4.4$\pm$0.0 \\
ELM-FBPINN & (40, 20) & (1, 12) & 2.9 & 3.5e+18 & LSQR & 1.4$\pm$1.7e-05 & 4.5$\pm$0.0 \\
ELM-FBPINN & (80, 40) & (1, 12) & 2.9 & 1.6e+19 & LSQR & 1.3$\pm$2.5e-04 & 5.1$\pm$0.0 \\
ELM-FBPINN & (160, 80) & (1, 12) & 2.9 & 1.4e+19 & LSQR & 2.6$\pm$2.0e-03 & 5.8$\pm$0.0 \\
ELM-FBPINN & (320, 160) & (1, 12) & 2.9 & 8.0e+19 & LSQR & 2.1$\pm$0.9e-02 & 8.0$\pm$0.1 \\
ELM-FBPINN & (20, [1, 5, 10]) & (1, 12) & 2.9 & 6.3e+15 & LSQR & 8.8$\pm$7.9e-07 & 6.5$\pm$0.1 \\
ELM-FBPINN & (40, [1, 10, 20]) & (1, 12) & 2.9 & 9.3e+14 & LSQR & 7.2$\pm$4.3e-06 & 7.3$\pm$0.1 \\
ELM-FBPINN & (80, [1, 20, 40]) & (1, 12) & 2.9 & 7.1e+14 & LSQR & 3.3$\pm$1.1e-05 & 6.7$\pm$0.0 \\
ELM-FBPINN & (160, [1, 40, 80]) & (1, 12) & 2.9 & 6.1e+14 & LSQR & 7.2$\pm$2.6e-05 & 8.6$\pm$0.0 \\
ELM-FBPINN & (320, [1, 80, 160]) & (1, 12) & 2.9 & 7.9e+14 & LSQR & 3.1$\pm$1.2e-04 & 1.6$\pm$0.0e+01 \\
\bottomrule
\end{tabular}

}
\captiontabscaling{$\omega_0$}{1D damped harmonic oscillator}{for}
\label{tab:ch3-harm-j-multiscale}
\end{table}

% ============================================================
\subsubsection{Summary of 1D findings}
The 1D experiments support four consistent conclusions:
(i) localisation is essential on this oscillatory problem (PINN fails; decomposed methods succeed),
(ii) for sufficient localisation and local width, ELM-FBPINN provides a faster route to accurate solutions than gradient-trained FBPINN,
(iii) increasing $J_1$ expands the trial space primarily via spatially distinct degrees of freedom, whereas increasing $K$ eventually introduces correlated features and ill-conditioning with diminishing \rev{reduction of the test error},
and (iv) using multilevel models with coarse levels further improves performance by aiding global subdomain communication.

\subsection{2D -- Multi-scale Laplacian problem}
\label{sec:ch3_num_res_laplace}

We next evaluate all models with a 2D multi-scale Laplacian problem, as defined in \cite{Dolean:MDD:2024}. Compared to the 1D oscillator, this problem combines higher dimensionality with multi-frequency spatial structure, increasing both representational and optimisation difficulty. 

\paragraph{Problem definition}

The problem is given by
\begin{equation}
\label{eq:laplace_problem}
\left\{
\begin{array}{ll}
- \Delta u = f, & \quad \text{in } \Omega = [0,1]^2, \\
u = 0, & \quad \text{on } \partial \Omega.
\end{array}
\right.
\end{equation}
with a source term given by
\begin{equation} \label{eq:multiscale_laplace_f}
f(x_1,x_2) = \frac{2}{n_\omega}\sum_{i=1}^{n_\omega} \big( \omega_i \pi \big)^2 \sin({\omega_i \pi x_1}) \, \sin({\omega_i \pi x_2}).
\end{equation}
In this case the exact solution is
\begin{equation}
u(x_1,x_2) = \frac{1}{n_\omega}\sum_{i=1}^{n_\omega} \sin({\omega_i \pi x_1}) \, \sin({\omega_i \pi x_2}).
\end{equation}
We choose $\omega_i$ such that the solution is highly multi-scale, with the number of multi-scale components controlled by the parameter $n_\omega$, and we vary $n_\omega$ and $\omega_i$ to control the problem difficulty.

\paragraph{Hard boundary constraints}
For this problem, we use a hard-constraining approach, as originally described by \cite{Lagaris1998}, to assert boundary conditions. Specifically, we use the following solution ansatz
\begin{equation}
\tilde{u}(x_1,x_2) = \tanh(\omega_{n_{\omega}} x_1)\tanh(\omega_{n_{\omega}}(1-x_1))\tanh(\omega_{n_{\omega}} x_2)\tanh(\omega_{n_{\omega}}(1-x_2))\,\hat{u}(x_1,x_2).
\end{equation}
It is straightforward to verify that this ansatz inherently satisfies the boundary condition in \Cref{eq:laplace_problem}. As a result, the boundary loss terms in \Cref{eq:multilevel_loss} can be omitted. Importantly, with this ansatz, \Cref{eq:multilevel_loss} remains quadratic in $\ba$, ensuring that the least-squares system in \Cref{eq:multilevel_ls_all_constraints} can still be formed. We also note that the $\tanh$ terms become sharper with increasing $\omega_{n_{\omega}}$, which is crucial to guarantee that the effect of the constraining operator is negligible at a distance of one solution half-cycle ($1/\omega_{n_{\omega}}$) from the boundary.

\paragraph{Physics-informed objective}
With the use of hard boundary constraints, the loss function reduces to
\begin{equation}
\mathcal{L}
=
\mathcal{L}_{\mathrm{phys}}=\frac{1}{N_I}
\sum_{i=1}^{N_I}
\big(
- \Delta \tilde{u} - f
\big)^2(\bx_i), \quad \bx_i \in \Omega.
\end{equation}

\subsubsection{Baseline}

\begin{figure}[h]
    \centering
    \includegraphics[width=\textwidth]{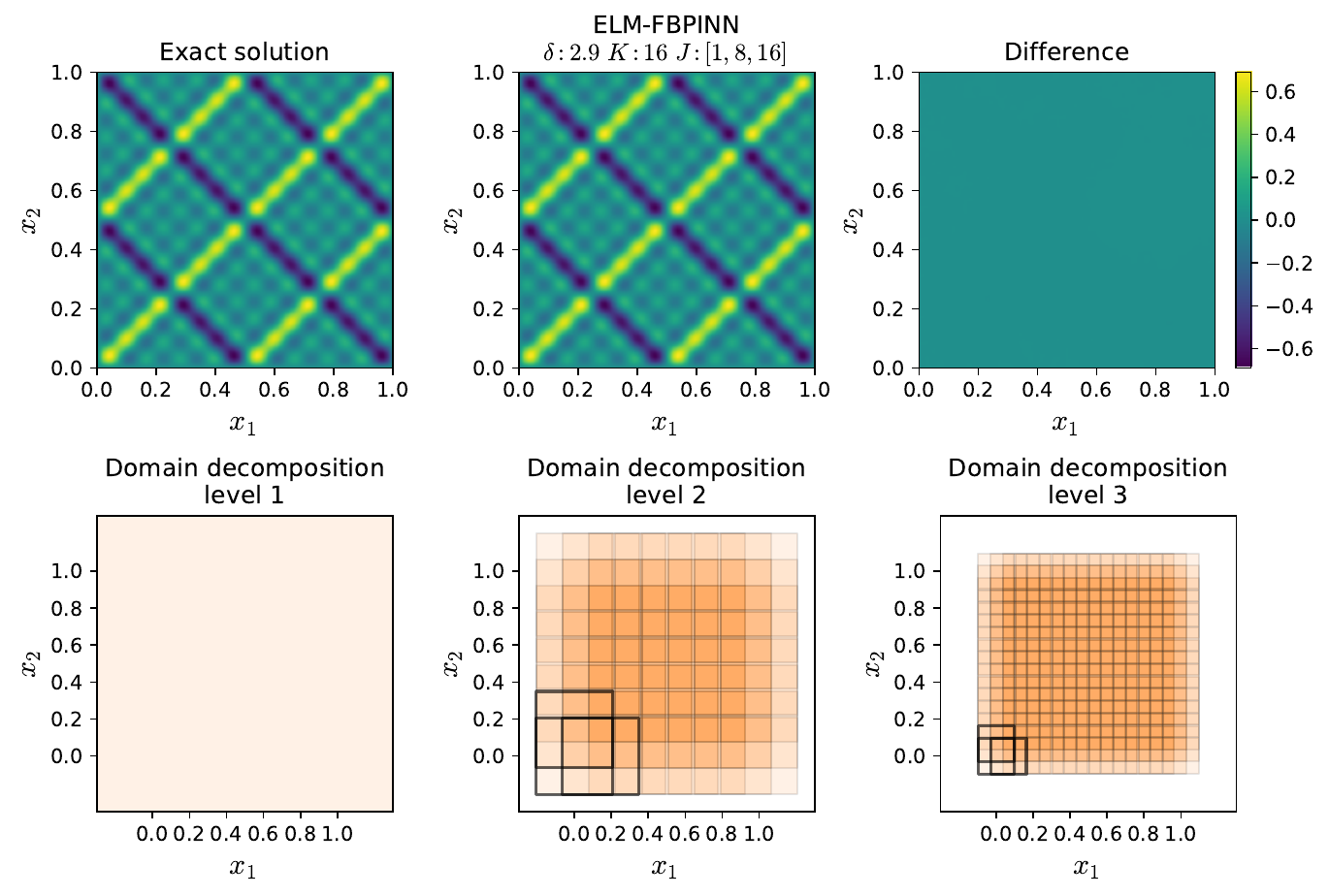}
    \caption{
    Model predictions versus the exact solution for the 2D multi-scale Laplacian problem with a multilevel ELM-FBPINN, where $J$ lists the number of subdomains along each input dimension for each level; multilevel domain decomposition.} 
    \label{fig:ch3-laplace-p0-preds}
\end{figure}

We establish a baseline configuration. We fix $n_\omega = 4$, choosing frequency components $\omega_1 = 4$, $\omega_2 = 8$, $\omega_3 = 12$, and $\omega_4 = 16$. The resulting solution (see \Cref{fig:ch3-laplace-p0-preds}) contains multiple interacting sine modes, introducing progressively finer spatial oscillations that are challenging for globally trained networks. The PINN uses $h=2$ hidden layers with $K=64$ basis functions per layer, matching the 1D configuration.  Both the FBPINN and ELM-FBPINN use a single hidden layer ($h=1$) with $K=16$, and a one-level decomposition with $J_1=16\times16=256$ subdomains with a width ratio $\delta=2.9$. {All models use a weight scaling of $R=1$ and} are trained on $80\times80=6400$ regularly spaced collocation points and tested on $M=96\times96=9216$ regularly spaced points.

\begin{figure}
    \centering
    \includegraphics[width=0.9\linewidth]{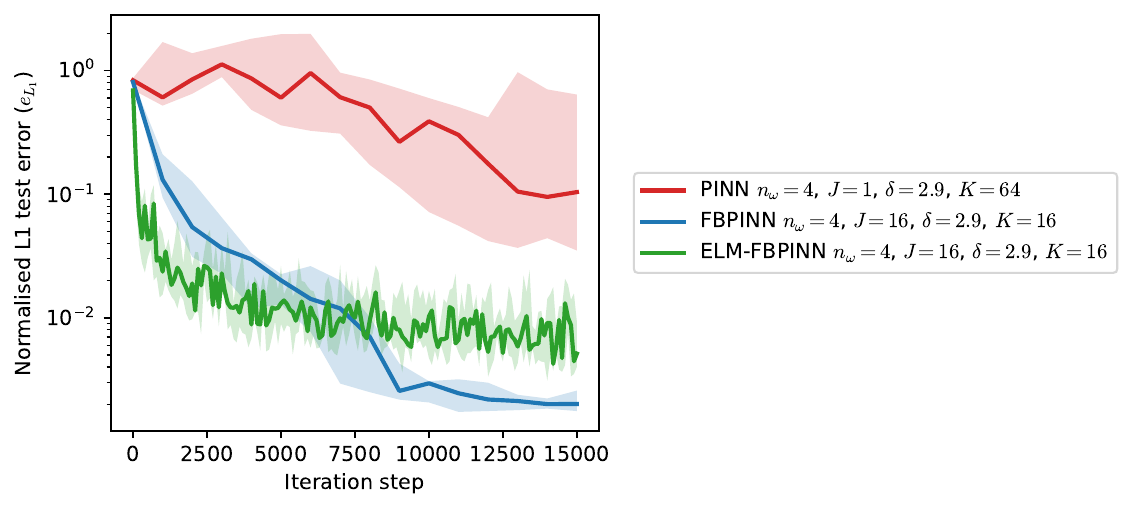}
\captionfigbase{2D multi-scale Laplacian}{along each input dimension for}
    \label{fig:ch3-laplace-p0-convergence}
\end{figure}

\Cref{fig:ch3-laplace-p0-convergence} shows the convergence curves. The PINN does not converge fully within the allotted budget, reflecting the difficulty of fitting multi-scale structure with a single global network trained by gradient descent. Both FBPINN and ELM-FBPINN converge to acceptable $e_{L^1}$ values, with FBPINN achieving the lowest final error. The improvement over PINN highlights the benefit of domain localisation. ELM-FBPINN exhibits rapid initial convergence, corresponding to the early iterations of the least-squares solver efficiently reducing dominant residual components, followed by a slower phase as finer-scale interactions dominate. \Cref{tab:ch3-laplace-p0-results} reports final $e_{L^1}$ values and wall-clock times. The PINN and FBPINN are nearly two orders of magnitude slower in runtime than the ELM-FBPINN, reflecting increased expense of nonlinear optimisation in the higher-dimensional multi-scale setting.

\begin{table}[h]
\setlength{\tabcolsep}{2pt}
\centering
\resizebox{\textwidth}{!}{
\begin{tabular}{cccccccc}
\toprule
Model & $(n_{\omega},J)$ & $(h,K)$ & $\delta$ & $\kappa(M)$ & Optimiser & $e_{L_1}$ & Time (s) \\
\midrule
PINN & (4, 1) & (2, 64) & 2.9 & N/A & Adam & 8.1$\pm$38.4e-02 & 1.8$\pm$0.0e+03 \\
FBPINN & (4, 16) & (1, 16) & 2.9 & N/A & Adam & 1.8$\pm$0.4e-03 & 2.4$\pm$0.0e+03 \\
ELM-FBPINN & (4, 16) & (1, 16) & 2.9 & 1.4e+09 & LSQR & 3.8$\pm$1.2e-03 & 6.2$\pm$0.0e+01 \\
\bottomrule
\end{tabular}

}
\captiontabbase{2D multi-scale Laplacian}{along each input dimension for}
\label{tab:ch3-laplace-p0-results}
\end{table}

We also show the model predictions in \Cref{fig:ch3-laplace-p0-convergence} for the same ELM-FBPINN, but with $L=3$ and  $J = [1, 8, 16]$ subdomains, where $J$ lists the number of subdomains along each input dimension for each level (i.e., $J_1=1\times1,J_2=8\times8$, and $J_3=16\times16$). The multilevel solution is visually indistinguishable from the exact solution, and recovers not only the dominant spatial oscillations but also the finer-scale components.

\subsubsection{Varying number of basis functions}

We vary the number of basis functions $K$ of the baseline $L=1$ model to probe how each method exploits increased local capacity while keeping domain decomposition and sampling fixed. This ablation separates (i) expressivity limitations at small $K$ from (ii) redundancy and conditioning effects at larger $K$.  \rev{As for the 1D problem, for the baseline PINN which uses $h=2$, varying $K$ also varies the number of basis functions per hidden layer as we use $K$ basis functions per hidden layer.}

\begin{figure}[t]
    \centering
    \includegraphics[width=\linewidth]{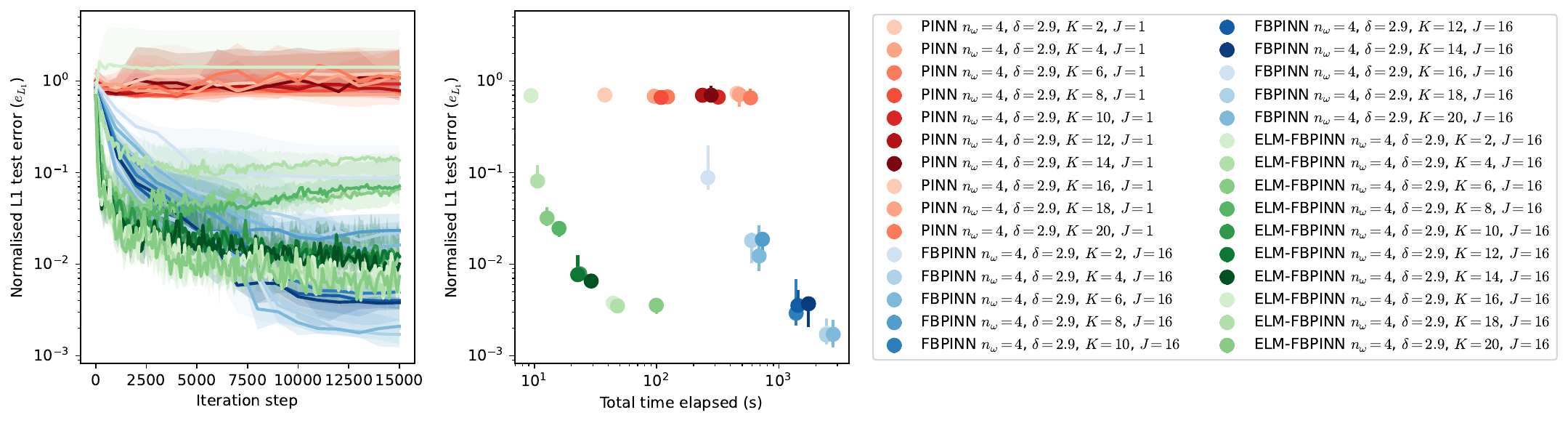}
\captionfigK{2D multi-scale Laplacian}{along each input dimension for}
    \label{fig:ch3-laplace-c-convergence}
\end{figure}

\paragraph{Key observations}
\Cref{fig:ch3-laplace-c-convergence} shows the convergence behaviour and \Cref{tab:ch3-laplace-C-results} summarises the final errors and runtimes. All PINN configurations exhibit similar non-convergent trends, indicating that increasing width alone does not resolve the optimisation difficulty associated with multi-scale structure in a single global network. Both decomposed methods converge for moderate $K$. FBPINN improves steadily as $K$ increases and achieves the lowest final errors, demonstrating that additional trainable local basis functions enhance expressivity in the 2D multi-scale setting. ELM-FBPINN also improves with increasing $K$, although very small values ($K=2,4$) lead to instability or large variance due to insufficient capacity and sensitivity to random feature alignment. ELM-FBPINN consistently exhibits faster initial convergence \rev{in terms of number of iteration steps} than FBPINN, reflecting the efficiency of solving a single least-squares problem for the dominant coefficients. For moderate to large $K$, its final errors remain within approximately one order of magnitude of FBPINN while requiring nearly two orders of magnitude less wall-clock time, highlighting a precision–efficiency trade-off.

\begin{table}[h]
\setlength{\tabcolsep}{2pt}
\centering
\resizebox{\textwidth}{!}{
\begin{tabular}{cccccccc}
\toprule
Model & $(n_{\omega},J)$ & $(h,K)$ & $\delta$ & $\kappa(M)$ & Optimiser & $e_{L_1}$ & Time (s) \\
\midrule
PINN & (4, 1) & (2, 2) & 2.9 & N/A & Adam & 7.0$\pm$0.3e-01 & 3.8$\pm$0.0e+01 \\
PINN & (4, 1) & (2, 4) & 2.9 & N/A & Adam & 6.8$\pm$0.9e-01 & 9.6$\pm$0.0e+01 \\
PINN & (4, 1) & (2, 6) & 2.9 & N/A & Adam & 6.7$\pm$0.5e-01 & 1.2$\pm$0.1e+02 \\
PINN & (4, 1) & (2, 8) & 2.9 & N/A & Adam & 6.6$\pm$0.8e-01 & 1.1$\pm$0.0e+02 \\
PINN & (4, 1) & (2, 10) & 2.9 & N/A & Adam & 6.7$\pm$1.2e-01 & 3.2$\pm$0.1e+02 \\
PINN & (4, 1) & (2, 12) & 2.9 & N/A & Adam & 7.0$\pm$1.1e-01 & 2.4$\pm$0.1e+02 \\
PINN & (4, 1) & (2, 14) & 2.9 & N/A & Adam & 7.0$\pm$2.0e-01 & 2.8$\pm$0.1e+02 \\
PINN & (4, 1) & (2, 16) & 2.9 & N/A & Adam & 7.4$\pm$1.5e-01 & 4.5$\pm$0.0e+02 \\
PINN & (4, 1) & (2, 18) & 2.9 & N/A & Adam & 7.1$\pm$2.1e-01 & 4.7$\pm$0.0e+02 \\
PINN & (4, 1) & (2, 20) & 2.9 & N/A & Adam & 6.6$\pm$2.4e-01 & 5.9$\pm$0.1e+02 \\
FBPINN & (4, 16) & (1, 2) & 2.9 & N/A & Adam & 8.8$\pm$13.2e-02 & 2.6$\pm$0.0e+02 \\
FBPINN & (4, 16) & (1, 4) & 2.9 & N/A & Adam & 1.8$\pm$1.2e-02 & 6.0$\pm$0.0e+02 \\
FBPINN & (4, 16) & (1, 6) & 2.9 & N/A & Adam & 1.2$\pm$1.8e-02 & 6.9$\pm$0.0e+02 \\
FBPINN & (4, 16) & (1, 8) & 2.9 & N/A & Adam & 1.9$\pm$1.1e-02 & 7.3$\pm$0.0e+02 \\
FBPINN & (4, 16) & (1, 10) & 2.9 & N/A & Adam & 2.9$\pm$4.7e-03 & 1.4$\pm$0.0e+03 \\
FBPINN & (4, 16) & (1, 12) & 2.9 & N/A & Adam & 3.5$\pm$2.5e-03 & 1.4$\pm$0.0e+03 \\
FBPINN & (4, 16) & (1, 14) & 2.9 & N/A & Adam & 3.7$\pm$1.8e-03 & 1.7$\pm$0.0e+03 \\
FBPINN & (4, 16) & (1, 16) & 2.9 & N/A & Adam & 1.8$\pm$0.4e-03 & 2.4$\pm$0.1e+03 \\
FBPINN & (4, 16) & (1, 18) & 2.9 & N/A & Adam & 1.7$\pm$1.2e-03 & 2.4$\pm$0.0e+03 \\
FBPINN & (4, 16) & (1, 20) & 2.9 & N/A & Adam & 1.7$\pm$1.3e-03 & 2.8$\pm$0.0e+03 \\
ELM-FBPINN & (4, 16) & (1, 2) & 2.9 & 1.9e+03 & LSQR & 6.9$\pm$0.0e-01 & 9.4$\pm$0.2 \\
ELM-FBPINN & (4, 16) & (1, 4) & 2.9 & 5.6e+04 & LSQR & 8.1$\pm$4.6e-02 & 1.1$\pm$0.0e+01 \\
ELM-FBPINN & (4, 16) & (1, 6) & 2.9 & 1.1e+06 & LSQR & 3.2$\pm$1.6e-02 & 1.3$\pm$0.0e+01 \\
ELM-FBPINN & (4, 16) & (1, 8) & 2.9 & 5.1e+17 & LSQR & 2.5$\pm$0.8e-02 & 1.6$\pm$0.0e+01 \\
ELM-FBPINN & (4, 16) & (1, 10) & 2.9 & 5.3e+07 & LSQR & 7.9$\pm$2.1e-03 & 2.4$\pm$0.0e+01 \\
ELM-FBPINN & (4, 16) & (1, 12) & 2.9 & 5.9e+08 & LSQR & 7.7$\pm$5.5e-03 & 2.3$\pm$0.0e+01 \\
ELM-FBPINN & (4, 16) & (1, 14) & 2.9 & 2.5e+10 & LSQR & 6.5$\pm$0.9e-03 & 2.9$\pm$0.1e+01 \\
ELM-FBPINN & (4, 16) & (1, 16) & 2.9 & 1.4e+09 & LSQR & 3.8$\pm$1.2e-03 & 4.4$\pm$0.0e+01 \\
ELM-FBPINN & (4, 16) & (1, 18) & 2.9 & 1.6e+10 & LSQR & 3.5$\pm$0.7e-03 & 4.8$\pm$0.1e+01 \\
ELM-FBPINN & (4, 16) & (1, 20) & 2.9 & 1.6e+11 & LSQR & 3.6$\pm$0.9e-03 & 9.9$\pm$0.4e+01 \\
\bottomrule
\end{tabular}

}
\captiontabK{2D multi-scale Laplacian}{along each input dimension for}
\label{tab:ch3-laplace-C-results}
\end{table}

\paragraph{Interpretation}
In the low-capacity regime, both decomposed methods are limited by insufficient local expressivity, with ELM-FBPINN particularly sensitive to random feature alignment. As $K$ increases, \rev{the test error} improves until diminishing returns appear. The persistent failure of PINN across all $K$ confirms that width alone does not compensate for optimisation difficulty in multi-scale, two-dimensional problems. Overall, FBPINN benefits from the flexibility of trainable hidden weights, whereas ELM-FBPINN offers competitive \rev{errors} at significantly reduced computational cost.

\subsubsection{Varying number of subdomains}

We vary the number of subdomains $J_1$ of the baseline $L=1$ model while fixing all other parameters for FBPINN and ELM-FBPINN. Here we use $J$ to denote the number of subdomains along each input dimension, i.e. $J_1=J\times J$.

\begin{figure}[h]
    \centering
    \includegraphics[width=\linewidth]{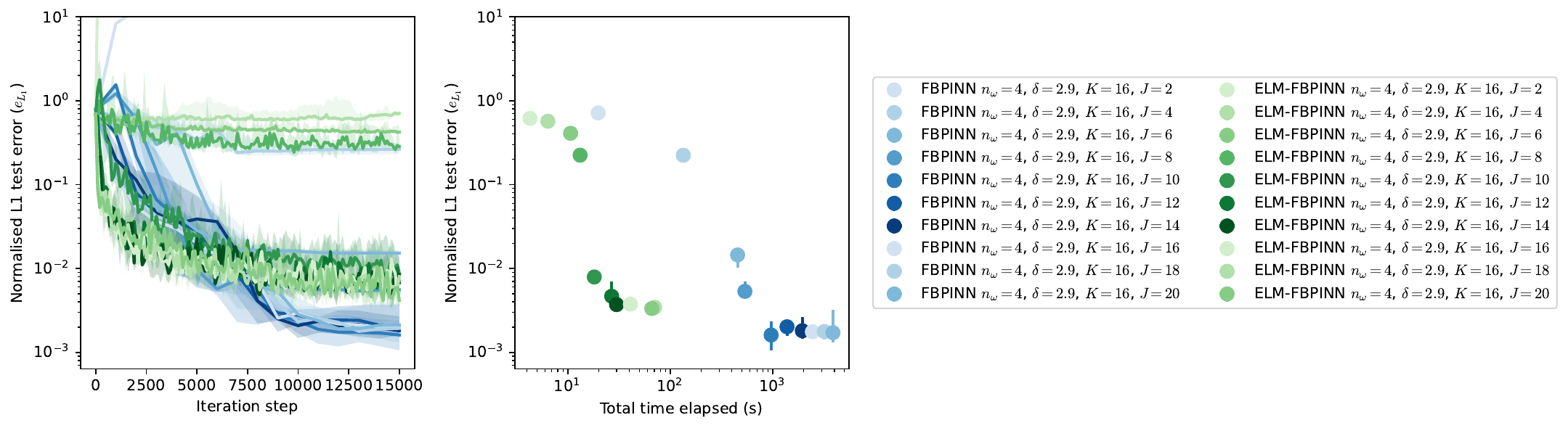}
\captionfigJ{2D multi-scale Laplacian}{along each input dimension for}
    \label{fig:ch3-laplace-j-convergence}
\end{figure}

\paragraph{Key observations}
\Cref{fig:ch3-laplace-j-convergence} shows the convergence behaviour. For small $J$, both methods struggle, with $J=2$ failing for both and $J=4$ exhibiting instability. In this regime, localisation is too coarse to isolate the multi-frequency structure of the solution. For moderate and large $J$, both methods converge reliably. FBPINN consistently attains lower $e_{L^1}$ values than ELM-FBPINN, with the gap remaining approximately one order of magnitude across most tested values. In contrast to the 1D case, increasing localisation alone is insufficient for ELM-FBPINN to match the \rev{errors} of fully trainable subdomain networks. ELM-FBPINN exhibits faster initial convergence \rev{in terms of number of iteration steps} for $J \gtrsim 10$, reflecting the efficiency of the global least-squares solve. The subsequent slower decay phase suggests residual inter-subdomain interactions that are more difficult to capture with fixed random features. Final errors and runtimes are summarised in \Cref{tab:ch3-laplace-J-results}. While ELM-FBPINN retains a significant runtime advantage of nearly two orders of magnitude, this gain comes at the cost of reduced precision.

\begin{table}[h!]
\setlength{\tabcolsep}{2pt}
\centering
\resizebox{\textwidth}{!}{
\begin{tabular}{cccccccc}
\toprule
Model & $(n_{\omega},J)$ & $(h,K)$ & $\delta$ & $\kappa(M)$ & Optimiser & $e_{L_1}$ & Time (s) \\
\midrule
FBPINN & (4, 2) & (1, 16) & 2.9 & N/A & Adam & 7.2$\pm$1.5e-01 & 2.0$\pm$0.0e+01 \\
FBPINN & (4, 4) & (1, 16) & 2.9 & N/A & Adam & 2.2$\pm$0.2e-01 & 1.3$\pm$0.0e+02 \\
FBPINN & (4, 6) & (1, 16) & 2.9 & N/A & Adam & 1.5$\pm$0.6e-02 & 4.5$\pm$0.2e+02 \\
FBPINN & (4, 8) & (1, 16) & 2.9 & N/A & Adam & 5.3$\pm$2.1e-03 & 5.3$\pm$0.1e+02 \\
FBPINN & (4, 10) & (1, 16) & 2.9 & N/A & Adam & 1.6$\pm$1.3e-03 & 9.6$\pm$0.2e+02 \\
FBPINN & (4, 12) & (1, 16) & 2.9 & N/A & Adam & 2.0$\pm$0.5e-03 & 1.4$\pm$0.0e+03 \\
FBPINN & (4, 14) & (1, 16) & 2.9 & N/A & Adam & 1.8$\pm$1.2e-03 & 1.9$\pm$0.0e+03 \\
FBPINN & (4, 16) & (1, 16) & 2.9 & N/A & Adam & 1.8$\pm$0.4e-03 & 2.4$\pm$0.0e+03 \\
FBPINN & (4, 18) & (1, 16) & 2.9 & N/A & Adam & 1.8$\pm$0.5e-03 & 3.2$\pm$0.1e+03 \\
FBPINN & (4, 20) & (1, 16) & 2.9 & N/A & Adam & 1.7$\pm$1.9e-03 & 3.9$\pm$0.0e+03 \\
ELM-FBPINN & (4, 2) & (1, 16) & 2.9 & 1.3e+10 & LSQR & 6.2$\pm$0.0e-01 & 4.3$\pm$0.0 \\
ELM-FBPINN & (4, 4) & (1, 16) & 2.9 & 7.6e+08 & LSQR & 5.7$\pm$0.1e-01 & 6.4$\pm$0.0 \\
ELM-FBPINN & (4, 6) & (1, 16) & 2.9 & 6.6e+08 & LSQR & 4.1$\pm$0.1e-01 & 1.1$\pm$0.0e+01 \\
ELM-FBPINN & (4, 8) & (1, 16) & 2.9 & 1.1e+09 & LSQR & 2.2$\pm$0.5e-01 & 1.3$\pm$0.0e+01 \\
ELM-FBPINN & (4, 10) & (1, 16) & 2.9 & 2.6e+09 & LSQR & 7.9$\pm$1.4e-03 & 1.8$\pm$0.0e+01 \\
ELM-FBPINN & (4, 12) & (1, 16) & 2.9 & 1.1e+09 & LSQR & 4.7$\pm$3.3e-03 & 2.7$\pm$0.0e+01 \\
ELM-FBPINN & (4, 14) & (1, 16) & 2.9 & 1.8e+09 & LSQR & 3.7$\pm$0.8e-03 & 3.0$\pm$0.0e+01 \\
ELM-FBPINN & (4, 16) & (1, 16) & 2.9 & 1.4e+09 & LSQR & 3.8$\pm$1.2e-03 & 4.1$\pm$0.0e+01 \\
ELM-FBPINN & (4, 18) & (1, 16) & 2.9 & 5.7e+09 & LSQR & 3.4$\pm$0.8e-03 & 7.0$\pm$0.0e+01 \\
ELM-FBPINN & (4, 20) & (1, 16) & 2.9 & 3.3e+09 & LSQR & 3.3$\pm$0.5e-03 & 6.6$\pm$0.0e+01 \\
\bottomrule
\end{tabular}

}
\captiontabJ{2D multi-scale Laplacian}{along each input dimension for}
\label{tab:ch3-laplace-J-results}
\end{table}

\paragraph{Interpretation}
In the 2D multi-scale setting, localisation alone is insufficient for fixed random features to match the flexibility of trainable hidden weights. While increasing $J$ improves \rev{the test loss errors} for both methods, the persistent performance gap indicates that adaptivity within subdomains becomes increasingly important as spatial complexity grows. ELM-FBPINN remains computationally attractive, but its fixed feature basis limits its ability to fully capture interacting frequency components in two dimensions.

\subsubsection{Scaling solution frequency with subdomains}

We investigate robustness under increasing spectral complexity by varying the number of frequency components $n_\omega$ in the exact solution. Specifically, our most complex problem is defined as $n_{\omega}=7$ with $\omega_i \in \{4, 8, 12, 16, 20, 24, 28\}, \,\, i = 1, \dots, 7$, with simpler problems defined by taking the first $n_{\omega}$ elements of this set. To maintain approximately constant resolution per subdomain, we scale the number of subdomains $J$ alongside $n_\omega$ (and the number of collocation points proportionally). This isolates the effect of increasing global complexity while preserving local resolution. In this section, we also assess whether including multiple levels can improve scaling performance.

\begin{figure}[t]
    \centering
    \includegraphics[width=\linewidth]{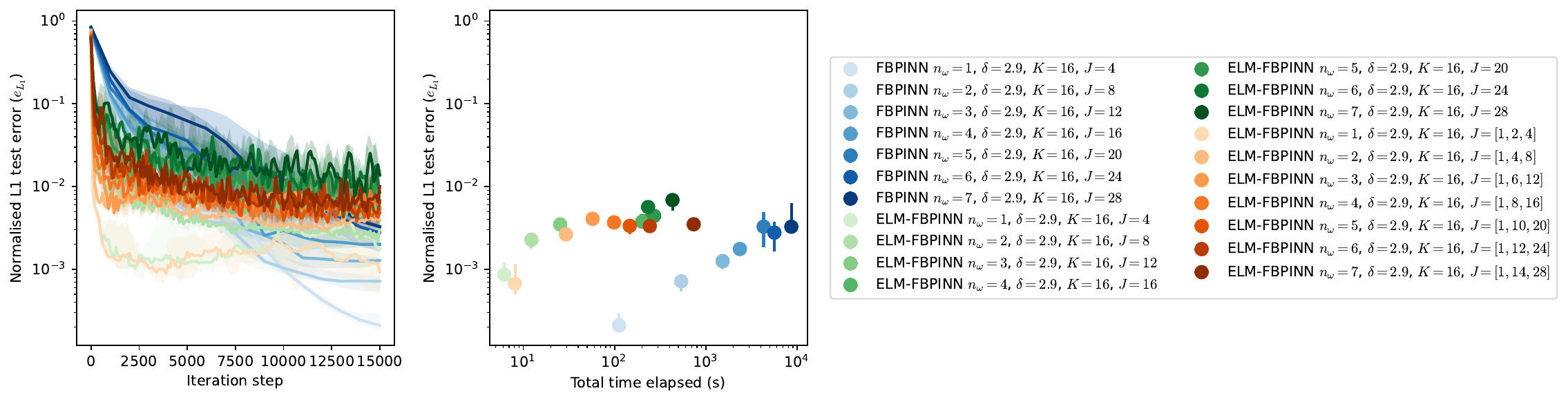}
\captionfigscaling{number of frequency components $n_\omega$}{2D multi-scale Laplacian}{along each input dimension for}
    \label{fig:ch3-laplace-scaling-convergence}
\end{figure}

\paragraph{Key observations}
\Cref{fig:ch3-laplace-scaling-convergence} shows the convergence behaviour and \Cref{tab:ch3-laplace-low_to_high-results} reports final errors. ELM-FBPINN and multilevel ELM-FBPINNs consistently exhibit the fastest initial convergence across all values of $n_\omega$, reflecting the efficiency of the global least-squares solve. Although convergence curves appear visually distinct, the absolute differences in final $e_{L^1}$ values across all models tested remain modest for most $n_\omega$, with the multilevel ELM-FBPINN most accurate on average, particularly for the highest frequencies tested. The max–min range of the ELM-FBPINN is slightly smaller than the FBPINN, indicating marginally more uniform performance as spectral richness increases. All methods maintain stable convergence as frequency increases, with no catastrophic degradation observed.

\begin{table}[h!]
\setlength{\tabcolsep}{2pt}
\centering
\resizebox{\textwidth}{!}{
\begin{tabular}{cccccccc}
\toprule
Model & $(n_{\omega},J)$ & $(h,K)$ & $\delta$ & $\kappa(M)$ & Optimiser & $e_{L_1}$ & Time (s) \\
\midrule
FBPINN & (1, 4) & (1, 16) & 2.9 & N/A & Adam & 2.1$\pm$1.1e-04 & 1.1$\pm$0.0e+02 \\
FBPINN & (2, 8) & (1, 16) & 2.9 & N/A & Adam & 7.2$\pm$2.7e-04 & 5.4$\pm$0.0e+02 \\
FBPINN & (3, 12) & (1, 16) & 2.9 & N/A & Adam & 1.3$\pm$0.3e-03 & 1.5$\pm$0.0e+03 \\
FBPINN & (4, 16) & (1, 16) & 2.9 & N/A & Adam & 1.8$\pm$0.4e-03 & 2.4$\pm$0.1e+03 \\
FBPINN & (5, 20) & (1, 16) & 2.9 & N/A & Adam & 3.3$\pm$3.1e-03 & 4.3$\pm$0.0e+03 \\
FBPINN & (6, 24) & (1, 16) & 2.9 & N/A & Adam & 2.8$\pm$2.1e-03 & 5.6$\pm$0.0e+03 \\
FBPINN & (7, 28) & (1, 16) & 2.9 & N/A & Adam & 3.3$\pm$3.3e-03 & 8.6$\pm$0.2e+03 \\
ELM-FBPINN & (1, 4) & (1, 16) & 2.9 & 1.0e+09 & LSQR & 8.6$\pm$5.0e-04 & 6.2$\pm$0.0 \\
ELM-FBPINN & (2, 8) & (1, 16) & 2.9 & 1.3e+09 & LSQR & 2.3$\pm$0.9e-03 & 1.2$\pm$0.0e+01 \\
ELM-FBPINN & (3, 12) & (1, 16) & 2.9 & 1.3e+09 & LSQR & 3.5$\pm$0.5e-03 & 2.5$\pm$0.0e+01 \\
ELM-FBPINN & (4, 16) & (1, 16) & 2.9 & 1.4e+09 & LSQR & 3.8$\pm$1.2e-03 & 2.0$\pm$0.1e+02 \\
ELM-FBPINN & (5, 20) & (1, 16) & 2.9 & 2.8e+09 & LSQR & 4.4$\pm$0.7e-03 & 2.7$\pm$0.5e+02 \\
ELM-FBPINN & (6, 24) & (1, 16) & 2.9 & 4.0e+09 & LSQR & 5.7$\pm$1.1e-03 & 2.3$\pm$0.1e+02 \\
ELM-FBPINN & (7, 28) & (1, 16) & 2.9 & 4.9e+09 & LSQR & 6.9$\pm$1.8e-03 & 4.3$\pm$0.2e+02 \\
ELM-FBPINN & (1, [1, 2, 4]) & (1, 16) & 2.9 & 3.8e+16 & LSQR & 6.7$\pm$6.7e-04 & 8.1$\pm$0.1 \\
ELM-FBPINN & (2, [1, 4, 8]) & (1, 16) & 2.9 & 3.0e+16 & LSQR & 2.6$\pm$0.4e-03 & 2.9$\pm$0.1e+01 \\
ELM-FBPINN & (3, [1, 6, 12]) & (1, 16) & 2.9 & 2.3e+16 & LSQR & 4.1$\pm$0.5e-03 & 5.8$\pm$0.1e+01 \\
ELM-FBPINN & (4, [1, 8, 16]) & (1, 16) & 2.9 & 2.2e+16 & LSQR & 3.7$\pm$1.5e-03 & 9.9$\pm$0.1e+01 \\
ELM-FBPINN & (5, [1, 10, 20]) & (1, 16) & 2.9 & 1.8e+16 & LSQR & 3.3$\pm$1.4e-03 & 1.5$\pm$0.0e+02 \\
ELM-FBPINN & (6, [1, 12, 24]) & (1, 16) & 2.9 & 2.1e+16 & LSQR & 3.3$\pm$1.1e-03 & 2.4$\pm$0.0e+02 \\
ELM-FBPINN & (7, [1, 14, 28]) & (1, 16) & 2.9 & 2.1e+16 & LSQR & 3.5$\pm$0.9e-03 & 7.4$\pm$0.6e+02 \\
\bottomrule
\end{tabular}

}
\captiontabscaling{$n_\omega$}{2D multi-scale Laplacian}{along each input dimension for}
\label{tab:ch3-laplace-low_to_high-results}
\end{table}

\paragraph{Interpretation}
All methods scale reasonably well as spectral richness increases, whilst the ELM-FBPINN and multilevel ELM-FBPINN remain nearly two orders of magnitude faster than the FBPINN; we note adaptive stopping for the least-squares solver could further reduce ELM-FBPINN cost without affecting accuracy.

\rev{
\subsection{2D -- Inhomogeneous Helmholtz equation}
\label{sec:ch3_num_res_helmholtz}}

\rev{We finally evaluate all models on a 2D Helmholtz problem with an inhomogeneous right-hand side on an L-shaped domain. Compared to the Laplacian problem, which has a simplifying tensor-product solution structure and homogeneous boundary conditions, we use a non-trivial manufactured solution which contains non-separable oscillatory solution components and spatially varying Dirichlet boundary conditions, yielding a more challenging and representative benchmark, particularly in the high-frequency regime.}

\rev{
\paragraph{Problem definition}
The problem is given by
\begin{equation}
\label{eq:helmholtz_problem}
\left\{
\begin{array}{ll}
\Delta u + k^2 u = f, & \quad \text{in } \Omega = [0,1]^2 \setminus [0.5,1]\times[0.5,1], \\
u = g, & \quad \text{on } \partial \Omega,
\end{array}
\right.
\end{equation}
where $k > 0$ is the wavenumber and is fixed at $k=1$. We define a manufactured solution
\begin{equation}
\label{eq:helmholtz_exact}
u(x_1,x_2)
=
\sin(\omega \pi x_1)\sin(2\omega \pi x_2)
+
\sin(3\omega \pi x_1 x_2),
\end{equation}
from which the source term is derived as
\begin{equation}
\label{eq:helmholtz_f}
\begin{aligned}
f(x_1,x_2)
=& \Delta u(x_1,x_2) + k^2 u(x_1,x_2) \\
=& -5(\omega \pi)^2
\sin(\omega \pi x_1)\sin(2\omega \pi x_2) -9(\omega \pi)^2(x_1^2+x_2^2)
\sin(3\omega \pi x_1 x_2) \\
&+k^2\Big(
\sin(\omega \pi x_1)\sin(2\omega \pi x_2)
+\sin(3\omega \pi x_1 x_2)
\Big),
\end{aligned}
\end{equation}
and non-homogeneous Dirichlet boundary conditions are imposed as $g = u|_{\partial\Omega}$. We choose this manufactured solution so that it contains both high-frequency and non-separable spatial components, yielding a non-trivial solution on the non-convex L-shaped domain with spatially varying boundary conditions. By increasing the frequency parameter $\omega$, we increase the oscillatory complexity of the solution and corresponding source term, thereby controlling the difficulty of the problem.
\paragraph{Physics-informed objective}
All models minimise the same residual-based loss,
\begin{equation}
\mathcal{L}
=
\mathcal{L}_{\mathrm{phys}}
+
\lambda_{\mathrm{bc}}
\mathcal{L}_{\mathrm{bc}},
\end{equation}
with
\begin{align}
\mathcal{L}_{\mathrm{phys}}
&=
\frac{1}{N_I}
\sum_{i=1}^{N_I}
\Big(
\Delta \hat u
+
k^2\hat{u}
-
f
\Big)^2(\bx_i), \quad \bx_i \in \Omega, \\
\mathcal{L}_{\mathrm{bc}}
&=
\frac{1}{N^{(\mathrm{bc})}_B}
\sum_{i=1}^{N^{(\mathrm{bc})}_B}
\big(
\hat{u}(\bx_i)-g(\bx_i)
\big)^2, \quad \bx_i \in \partial \Omega,
\end{align}
where $N_I$ and $N^{(\mathrm{bc})}_B$ denote the number of interior and boundary collocation points, respectively. We choose $\lambda_{\mathrm{bc}} = 25\omega^4$, which normalises all loss terms to the same physical unit and we find provides a reasonable balance between them as the solution frequency increases.}

\begin{figure}
    \centering
    \includegraphics[width=0.9\linewidth]{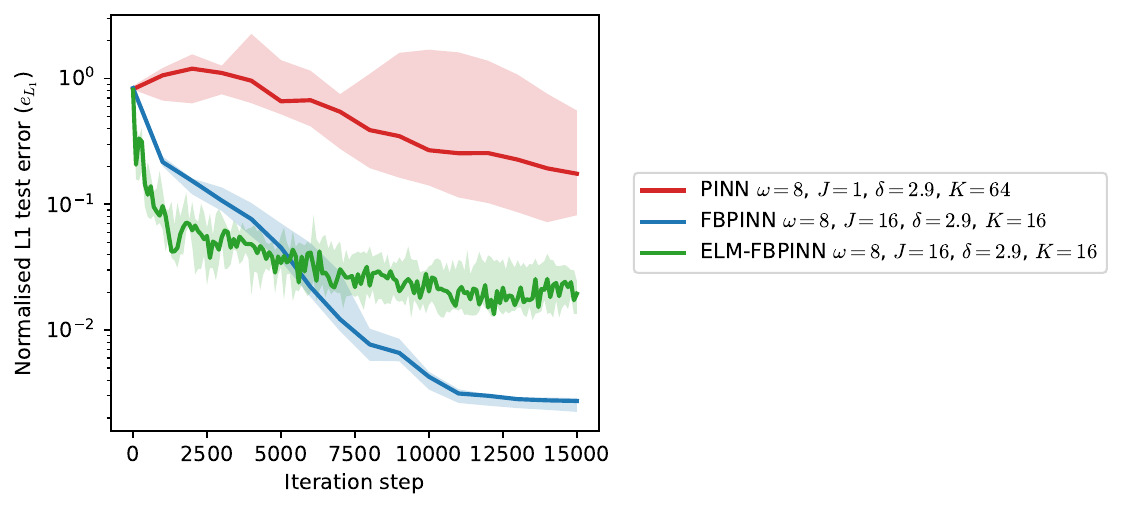}
\captionfigbase{2D inhomogenous Helmholtz equation}{along each input dimension for}
    \label{fig:ch3-helmholtz-p0-convergence}
\end{figure}

\begin{figure}[h]
    \centering
    \includegraphics[width=\textwidth]{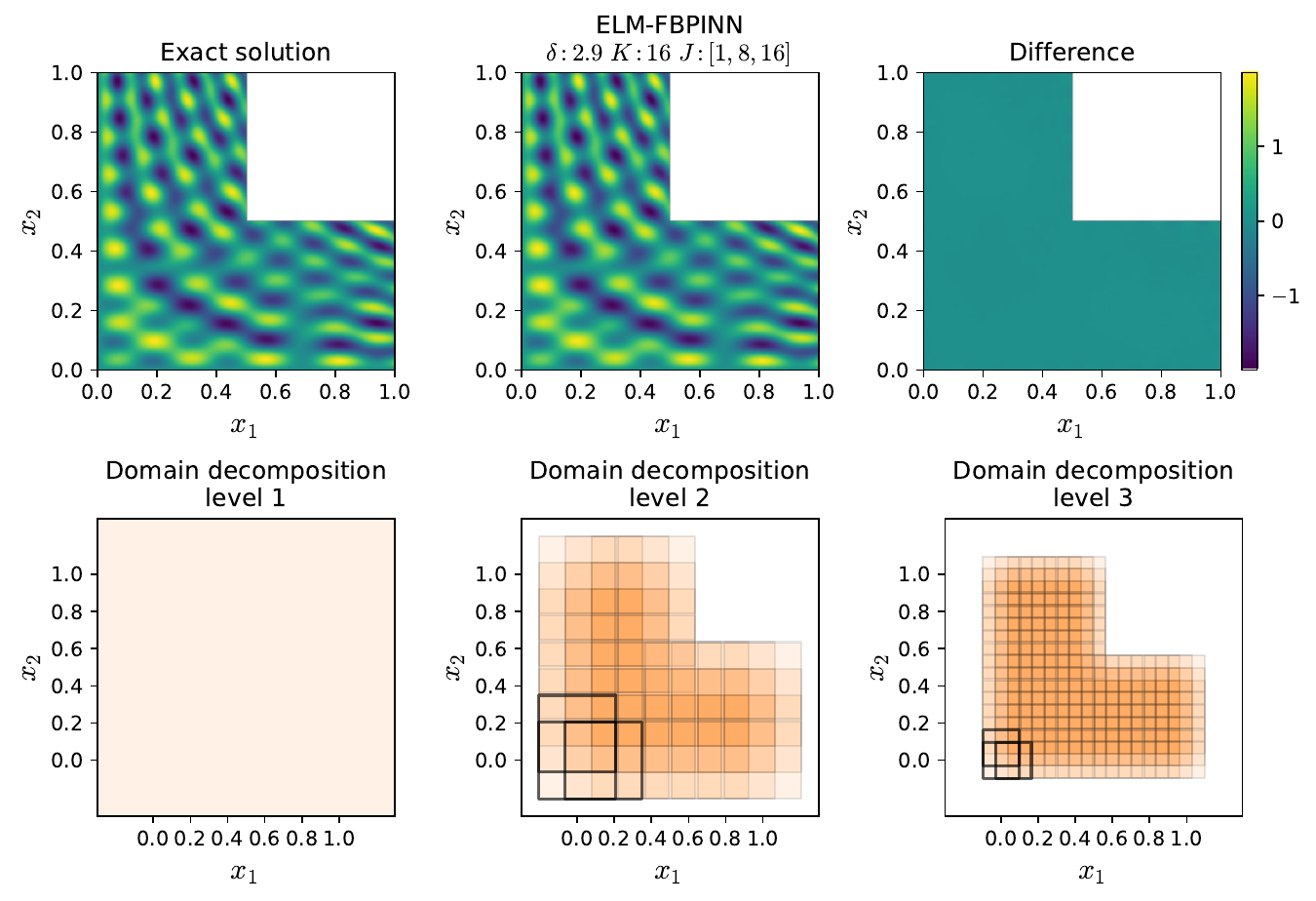}
    \caption{
    Model predictions versus the exact solution for the 2D inhomogenous Helmholtz equation problem with a multilevel ELM-FBPINN, where $J$ lists the number of subdomains along each input dimension for each level; multilevel domain decomposition.} 
    \label{fig:ch3-helmholtz-p0-preds}
\end{figure}

\rev{\subsubsection{Baseline}
We establish a baseline configuration. We fix $\omega = 8$, with the resulting solution shown in \cref{fig:ch3-helmholtz-p0-preds}. The PINN uses $h=2$ hidden layers with $K=64$ basis functions per layer. Both the FBPINN and ELM-FBPINN use $h=1$ with $K=16$, and a one-level rectangular decomposition with $J_1=$ $16\times16 - 8\times 8=192$ subdomains where $ 8\times 8$ subdomains outside of the L-shaped domain are discarded, and width ratio $\delta=2.9$. All models use a weight scaling of $R=1$. Interior collocation points are obtained from a regular $80\times80$ Cartesian grid over $[0,1]^2$, discarding $40\times40$ points outside the L-shaped domain, while boundary collocation points are sampled uniformly along each boundary segment (80 points on long segments and 40 on short segments). Test points are constructed analogously using a regular $96\times96$ Cartesian grid over $[0,1]^2$, again discarding points outside the L-shaped domain.\\
\Cref{fig:ch3-helmholtz-p0-convergence} shows convergence curves. As in the previous settings, the PINN exhibits slow and unstable convergence, failing to reach a competitive error within the allotted iteration budget. Both FBPINN and ELM-FBPINN converge reliably, with the FBPINN achieving a lower final error. The gap in convergence speed is particularly pronounced in the early iterations. Although the final error favours FBPINN, the time to solution is orders of magnitude lower for the ELM version (see Table \ref{tab:ch3-helmholtz-p0-results}).}

\begin{table}[h]
\setlength{\tabcolsep}{2pt}
\centering
\resizebox{\textwidth}{!}{
 \begin{tabular}{cccccccc}
\toprule
Model & $(\omega,J)$ & $(h,K)$ & $\delta$ & $\kappa(M)$ & Optimiser & $e_{L_1}$ & Time (s) \\
\midrule
PINN & (8, 1) & (2, 64) & 2.9 & N/A & Adam & 1.8$\pm$4.8e-01 & 1.1$\pm$0.0e+03 \\
FBPINN & (8, 16) & (1, 16) & 2.9 & N/A & Adam & 2.7$\pm$0.6e-03 & 1.5$\pm$0.0e+03 \\
ELM-FBPINN & (8, 16) & (1, 16) & 2.9 & 3.7e+08 & LSQR & 1.3$\pm$0.8e-02 & 3.7$\pm$0.0e+01 \\
\bottomrule
\end{tabular}
}
\captiontabbase{2D inhomogenous Helmholtz equation}{along each input dimension for}
\label{tab:ch3-helmholtz-p0-results}
\end{table}
\rev{We also show model predictions in \Cref{fig:ch3-helmholtz-p0-preds} for the same ELM-FBPINN, but with $L=3$ and  $J = [1, 8, 16]$ subdomains, where $J$ lists the number of subdomains along each input dimension for each level (i.e., $J_1=1\times1,J_2=8\times8-4\times4$, and $J_3=16\times16-8\times8$). The multilevel solution is visually indistinguishable from the exact solution.}

\rev{\subsubsection{Varying number of basis functions}
We vary the number of basis functions $K$ of the baseline $L=1$ model while keeping all other parameters fixed.  As before, since $h=2$ for the baseline PINN, varying $K$ alters the number of basis functions per hidden layer.} 

\begin{figure}[t]
    \centering
    \includegraphics[width=\linewidth]{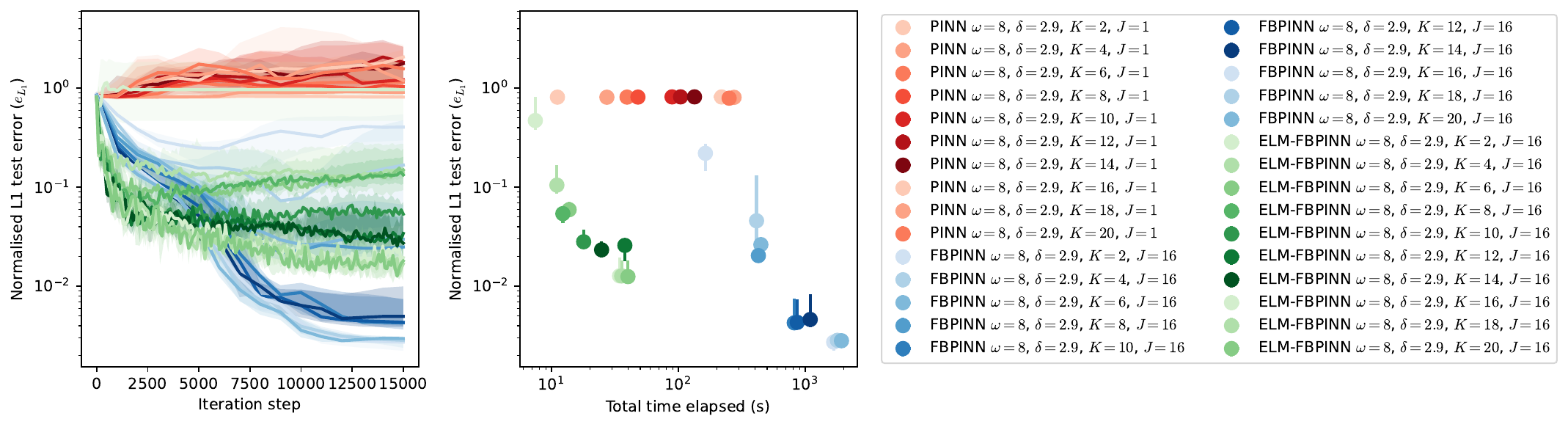}
\captionfigK{2D inhomogenous Helmholtz equation}{along each input dimension for}
    \label{fig:ch3-helmholtz-c-convergence}
\end{figure}

\rev{
\paragraph{Key observations}
\Cref{fig:ch3-helmholtz-c-convergence} shows the convergence behaviour and \Cref{tab:ch3-helmholtz-C-results} summarises the final errors and runtimes. Across all models, increasing $K$ improves approximation accuracy, although the effect is markedly different between methods. The PINN shows limited sensitivity to $K$, with only marginal improvements in the final error, suggesting that optimisation remains the dominant bottleneck. In contrast, both ELM-FBPINN and FBPINN benefit from increased basis size, with the latter exhibiting the most consistent and rapid error reduction. }

\begin{table}[h]
\setlength{\tabcolsep}{2pt}
\centering
\resizebox{\textwidth}{!}{
 \begin{tabular}{cccccccc}
\toprule
Model & $(\omega,J)$ & $(h,K)$ & $\delta$ & $\kappa(M)$ & Optimiser & $e_{L_1}$ & Time (s) \\
\midrule
PINN & (8, 1) & (2, 2) & 2.9 & N/A & Adam & 8.1$\pm$0.0e-01 & 1.1$\pm$0.0e+01 \\
PINN & (8, 1) & (2, 4) & 2.9 & N/A & Adam & 8.1$\pm$0.0e-01 & 2.7$\pm$0.0e+01 \\
PINN & (8, 1) & (2, 6) & 2.9 & N/A & Adam & 8.1$\pm$0.1e-01 & 3.9$\pm$0.0e+01 \\
PINN & (8, 1) & (2, 8) & 2.9 & N/A & Adam & 8.1$\pm$0.2e-01 & 4.8$\pm$0.2e+01 \\
PINN & (8, 1) & (2, 10) & 2.9 & N/A & Adam & 8.2$\pm$0.5e-01 & 8.9$\pm$0.1e+01 \\
PINN & (8, 1) & (2, 12) & 2.9 & N/A & Adam & 8.2$\pm$0.3e-01 & 1.0$\pm$0.0e+02 \\
PINN & (8, 1) & (2, 14) & 2.9 & N/A & Adam & 8.2$\pm$0.5e-01 & 1.3$\pm$0.0e+02 \\
PINN & (8, 1) & (2, 16) & 2.9 & N/A & Adam & 8.2$\pm$0.7e-01 & 2.2$\pm$0.2e+02 \\
PINN & (8, 1) & (2, 18) & 2.9 & N/A & Adam & 8.1$\pm$0.6e-01 & 2.7$\pm$0.2e+02 \\
PINN & (8, 1) & (2, 20) & 2.9 & N/A & Adam & 7.9$\pm$0.8e-01 & 2.5$\pm$0.0e+02 \\
FBPINN & (8, 16) & (1, 2) & 2.9 & N/A & Adam & 2.2$\pm$1.3e-01 & 1.6$\pm$0.0e+02 \\
FBPINN & (8, 16) & (1, 4) & 2.9 & N/A & Adam & 4.6$\pm$10.4e-02 & 4.1$\pm$0.1e+02 \\
FBPINN & (8, 16) & (1, 6) & 2.9 & N/A & Adam & 2.6$\pm$1.0e-02 & 4.5$\pm$0.0e+02 \\
FBPINN & (8, 16) & (1, 8) & 2.9 & N/A & Adam & 2.0$\pm$0.8e-02 & 4.3$\pm$0.0e+02 \\
FBPINN & (8, 16) & (1, 10) & 2.9 & N/A & Adam & 4.3$\pm$3.5e-03 & 8.2$\pm$0.0e+02 \\
FBPINN & (8, 16) & (1, 12) & 2.9 & N/A & Adam & 4.3$\pm$3.7e-03 & 8.6$\pm$0.0e+02 \\
FBPINN & (8, 16) & (1, 14) & 2.9 & N/A & Adam & 4.6$\pm$4.4e-03 & 1.1$\pm$0.0e+03 \\
FBPINN & (8, 16) & (1, 16) & 2.9 & N/A & Adam & 2.7$\pm$0.6e-03 & 1.7$\pm$0.0e+03 \\
FBPINN & (8, 16) & (1, 18) & 2.9 & N/A & Adam & 2.8$\pm$0.5e-03 & 1.8$\pm$0.0e+03 \\
FBPINN & (8, 16) & (1, 20) & 2.9 & N/A & Adam & 2.8$\pm$0.4e-03 & 1.9$\pm$0.0e+03 \\
ELM-FBPINN & (8, 16) & (1, 2) & 2.9 & 1.0e+03 & LSQR & 4.7$\pm$4.3e-01 & 7.4$\pm$0.0 \\
ELM-FBPINN & (8, 16) & (1, 4) & 2.9 & 5.1e+04 & LSQR & 1.1$\pm$0.8e-01 & 1.1$\pm$0.0e+01 \\
ELM-FBPINN & (8, 16) & (1, 6) & 2.9 & 4.7e+05 & LSQR & 5.9$\pm$1.6e-02 & 1.4$\pm$0.0e+01 \\
ELM-FBPINN & (8, 16) & (1, 8) & 2.9 & 1.8e+07 & LSQR & 5.4$\pm$1.9e-02 & 1.2$\pm$0.0e+01 \\
ELM-FBPINN & (8, 16) & (1, 10) & 2.9 & 1.7e+07 & LSQR & 2.8$\pm$1.3e-02 & 1.8$\pm$0.0e+01 \\
ELM-FBPINN & (8, 16) & (1, 12) & 2.9 & 1.6e+08 & LSQR & 2.6$\pm$0.9e-02 & 3.8$\pm$0.1e+01 \\
ELM-FBPINN & (8, 16) & (1, 14) & 2.9 & 2.6e+09 & LSQR & 2.3$\pm$0.9e-02 & 2.5$\pm$0.0e+01 \\
ELM-FBPINN & (8, 16) & (1, 16) & 2.9 & 3.7e+08 & LSQR & 1.3$\pm$0.8e-02 & 3.4$\pm$0.0e+01 \\
ELM-FBPINN & (8, 16) & (1, 18) & 2.9 & 2.1e+09 & LSQR & 1.3$\pm$0.7e-02 & 3.6$\pm$0.0e+01 \\
ELM-FBPINN & (8, 16) & (1, 20) & 2.9 & 2.2e+10 & LSQR & 1.3$\pm$0.7e-02 & 4.0$\pm$0.0e+01 \\
\bottomrule
\end{tabular}

 }
\captiontabK{2D inhomogenous Helmholtz equation}{along each input dimension for}
\label{tab:ch3-helmholtz-C-results}
\end{table}

\begin{figure}[h!]
    \centering
    \includegraphics[width=\linewidth]{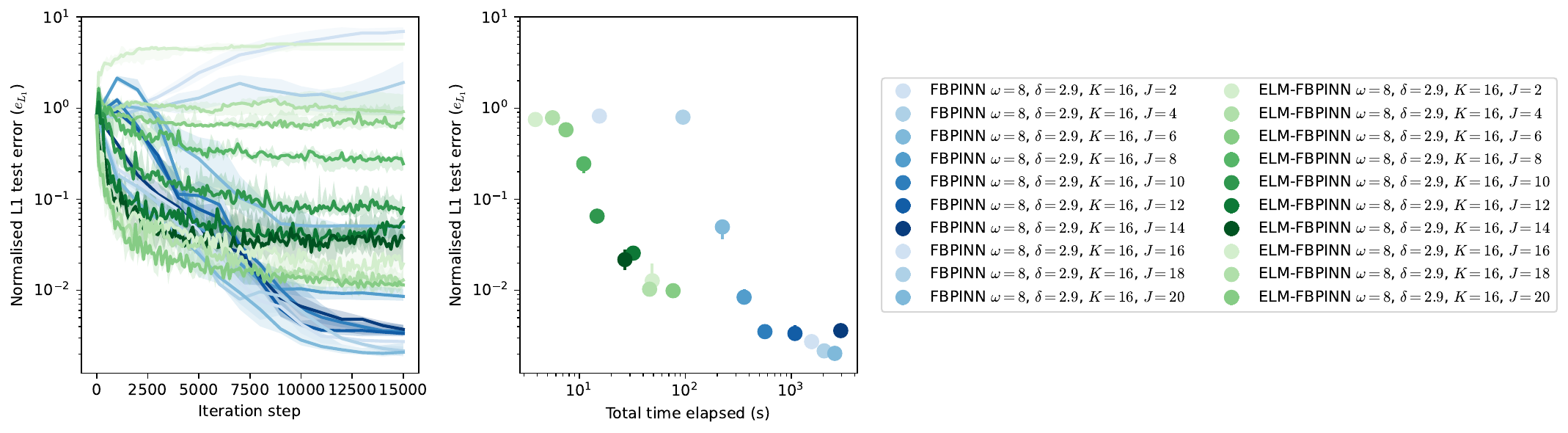}
\captionfigJ{2D inhomogenous Helmholtz equation}{along each input dimension for}
    \label{fig:ch3-helmholtz-j-convergence}
\end{figure}

\rev{\subsubsection{Varying number of subdomains}
We vary the number of subdomains $J_1$ of the baseline $L=1$ model while fixing all other parameters for FBPINN and ELM-FBPINN. Here we use $J$ to denote the number of subdomains along each input dimension, i.e. $J_1=J\times J-(J/2\times J/2)$.
\paragraph{Key observations}
\Cref{fig:ch3-helmholtz-j-convergence} shows the convergence behaviour and final errors and runtimes are summarised in \Cref{tab:ch3-helmholtz-J-results}. Increasing the number of subdomains leads to a clear improvement in accuracy for both FBPINN and ELM-FBPINN.  For small values of $J$, both methods struggle to resolve the oscillatory structure of the solution, while beyond a moderate threshold, the FBPINN consistently achieves lower errors than the ELM-FBPINN, but the latter being order of magnitude faster.}

\begin{table}[h!]
\setlength{\tabcolsep}{2pt}
\centering
 \resizebox{\textwidth}{!}{
 \begin{tabular}{cccccccc}
\toprule
Model & $(\omega,J)$ & $(h,K)$ & $\delta$ & $\kappa(M)$ & Optimiser & $e_{L_1}$ & Time (s) \\
\midrule
FBPINN & (8, 2) & (1, 16) & 2.9 & N/A & Adam & 8.2$\pm$0.5e-01 & 1.5$\pm$0.0e+01 \\
FBPINN & (8, 4) & (1, 16) & 2.9 & N/A & Adam & 8.0$\pm$0.4e-01 & 9.5$\pm$0.0e+01 \\
FBPINN & (8, 6) & (1, 16) & 2.9 & N/A & Adam & 5.0$\pm$1.6e-02 & 2.2$\pm$0.1e+02 \\
FBPINN & (8, 8) & (1, 16) & 2.9 & N/A & Adam & 8.4$\pm$2.6e-03 & 3.6$\pm$0.1e+02 \\
FBPINN & (8, 10) & (1, 16) & 2.9 & N/A & Adam & 3.5$\pm$0.3e-03 & 5.6$\pm$0.2e+02 \\
FBPINN & (8, 12) & (1, 16) & 2.9 & N/A & Adam & 3.4$\pm$1.1e-03 & 1.1$\pm$0.0e+03 \\
FBPINN & (8, 14) & (1, 16) & 2.9 & N/A & Adam & 3.6$\pm$1.0e-03 & 2.9$\pm$0.1e+03 \\
FBPINN & (8, 16) & (1, 16) & 2.9 & N/A & Adam & 2.7$\pm$0.6e-03 & 1.6$\pm$0.0e+03 \\
FBPINN & (8, 18) & (1, 16) & 2.9 & N/A & Adam & 2.2$\pm$0.5e-03 & 2.0$\pm$0.1e+03 \\
FBPINN & (8, 20) & (1, 16) & 2.9 & N/A & Adam & 2.0$\pm$0.4e-03 & 2.6$\pm$0.1e+03 \\
ELM-FBPINN & (8, 2) & (1, 16) & 2.9 & 7.5e+09 & LSQR & 7.5$\pm$0.0e-01 & 3.8$\pm$0.0 \\
ELM-FBPINN & (8, 4) & (1, 16) & 2.9 & 3.5e+08 & LSQR & 7.8$\pm$0.0e-01 & 5.6$\pm$0.0 \\
ELM-FBPINN & (8, 6) & (1, 16) & 2.9 & 3.7e+08 & LSQR & 5.8$\pm$0.7e-01 & 7.5$\pm$0.1 \\
ELM-FBPINN & (8, 8) & (1, 16) & 2.9 & 2.3e+08 & LSQR & 2.4$\pm$0.7e-01 & 1.1$\pm$0.0e+01 \\
ELM-FBPINN & (8, 10) & (1, 16) & 2.9 & 3.9e+08 & LSQR & 6.5$\pm$0.5e-02 & 1.5$\pm$0.0e+01 \\
ELM-FBPINN & (8, 12) & (1, 16) & 2.9 & 3.7e+08 & LSQR & 2.6$\pm$0.6e-02 & 3.2$\pm$0.0e+01 \\
ELM-FBPINN & (8, 14) & (1, 16) & 2.9 & 3.0e+08 & LSQR & 2.2$\pm$1.1e-02 & 2.7$\pm$0.0e+01 \\
ELM-FBPINN & (8, 16) & (1, 16) & 2.9 & 3.7e+08 & LSQR & 1.3$\pm$0.8e-02 & 4.9$\pm$0.2e+01 \\
ELM-FBPINN & (8, 18) & (1, 16) & 2.9 & 9.8e+08 & LSQR & 1.0$\pm$0.3e-02 & 4.6$\pm$0.1e+01 \\
ELM-FBPINN & (8, 20) & (1, 16) & 2.9 & 9.4e+08 & LSQR & 9.9$\pm$2.5e-03 & 7.7$\pm$0.1e+01 \\
\bottomrule
\end{tabular}

 }
\captiontabJ{2D inhomogenous Helmholtz equation}{along each input dimension for}
\label{tab:ch3-helmholtz-J-results}
\end{table}

\rev{\paragraph{Interpretation}
These results confirm the importance of localisation for oscillatory problems. By reducing the size of each subdomain, the local approximation task becomes simpler, allowing both models to better capture high-frequency behaviour. However, the gains diminish at larger $J$, indicating that beyond a certain point, further subdivision yields limited additional benefit.}

\begin{figure}[t]
    \centering
    \includegraphics[width=\linewidth]{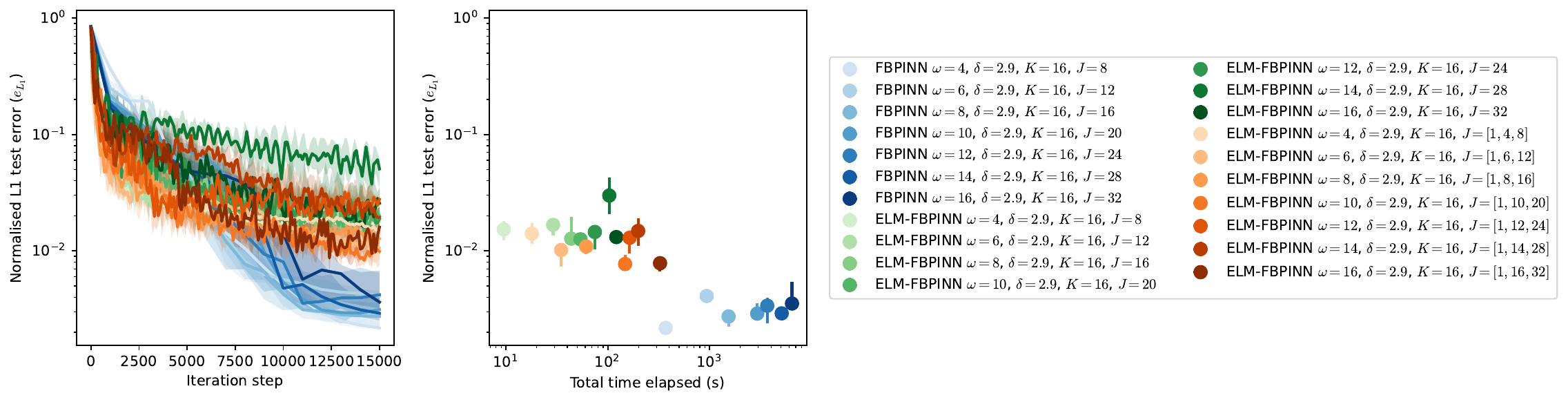}
\captionfigscaling{solution frequency $\omega$}{2D inhomogenous Helmholtz equation}{along each input dimension for}
    \label{fig:ch3-helmholtz-scaling-convergence}
\end{figure}

\rev{\subsubsection{Scaling solution frequency with subdomains}
We test scalability to harder (higher-frequency) solutions by increasing $\omega$ while simultaneously increasing $J_1$ (and the number of collocation points proportionally) to maintain approximately constant resolution per subdomain, for the baseline $L=1$ model. In this section, we also assess whether including multiple levels can improve scaling performance.
\paragraph{Key observations}
\Cref{fig:ch3-helmholtz-scaling-convergence} shows the convergence behaviour and \Cref{tab:ch3-helmholtz-scaling-results} reports final errors. As the solution frequency increases, all models exhibit slower convergence, in terms of iteration numbers and CPU time, and higher final errors.  The degradation in performance is gradual, indicating that the method scales reasonably well when the resolution per subdomain is preserved. Multilevel configurations further improve robustness, particularly at higher frequencies.}

\begin{table}[h!]
\setlength{\tabcolsep}{2pt}
\centering
 \resizebox{\textwidth}{!}{
 \begin{tabular}{cccccccc}
\toprule
Model & $(\omega,J)$ & $(h,K)$ & $\delta$ & $\kappa(M)$ & Optimiser & $e_{L_1}$ & Time (s) \\
\midrule
FBPINN & (4, 8) & (1, 16) & 2.9 & N/A & Adam & 2.2$\pm$0.3e-03 & 3.7$\pm$0.2e+02 \\
FBPINN & (6, 12) & (1, 16) & 2.9 & N/A & Adam & 4.1$\pm$0.8e-03 & 9.3$\pm$0.0e+02 \\
FBPINN & (8, 16) & (1, 16) & 2.9 & N/A & Adam & 2.7$\pm$0.6e-03 & 1.5$\pm$0.0e+03 \\
FBPINN & (10, 20) & (1, 16) & 2.9 & N/A & Adam & 2.9$\pm$0.8e-03 & 2.9$\pm$0.1e+03 \\
FBPINN & (12, 24) & (1, 16) & 2.9 & N/A & Adam & 3.4$\pm$1.6e-03 & 3.7$\pm$0.1e+03 \\
FBPINN & (14, 28) & (1, 16) & 2.9 & N/A & Adam & 2.9$\pm$0.6e-03 & 5.1$\pm$0.1e+03 \\
FBPINN & (16, 32) & (1, 16) & 2.9 & N/A & Adam & 3.5$\pm$2.2e-03 & 6.4$\pm$0.0e+03 \\
ELM-FBPINN & (4, 8) & (1, 16) & 2.9 & 1.8e+08 & LSQR & 1.5$\pm$0.6e-02 & 9.5$\pm$0.0 \\
ELM-FBPINN & (6, 12) & (1, 16) & 2.9 & 3.8e+08 & LSQR & 1.7$\pm$0.5e-02 & 2.9$\pm$0.0e+01 \\
ELM-FBPINN & (8, 16) & (1, 16) & 2.9 & 3.7e+08 & LSQR & 1.3$\pm$0.8e-02 & 4.4$\pm$0.3e+01 \\
ELM-FBPINN & (10, 20) & (1, 16) & 2.9 & 9.0e+08 & LSQR & 1.3$\pm$0.2e-02 & 5.4$\pm$0.0e+01 \\
ELM-FBPINN & (12, 24) & (1, 16) & 2.9 & 1.3e+09 & LSQR & 1.5$\pm$0.5e-02 & 7.5$\pm$0.0e+01 \\
ELM-FBPINN & (14, 28) & (1, 16) & 2.9 & 6.2e+08 & LSQR & 3.0$\pm$2.2e-02 & 1.0$\pm$0.0e+02 \\
ELM-FBPINN & (16, 32) & (1, 16) & 2.9 & 2.5e+09 & LSQR & 1.3$\pm$0.4e-02 & 1.2$\pm$0.0e+02 \\
ELM-FBPINN & (4, [1, 4, 8]) & (1, 16) & 2.9 & 3.5e+16 & LSQR & 1.4$\pm$0.6e-02 & 1.8$\pm$0.0e+01 \\
ELM-FBPINN & (6, [1, 6, 12]) & (1, 16) & 2.9 & 3.9e+16 & LSQR & 1.0$\pm$0.3e-02 & 3.5$\pm$0.0e+01 \\
ELM-FBPINN & (8, [1, 8, 16]) & (1, 16) & 2.9 & 5.3e+16 & LSQR & 1.1$\pm$0.4e-02 & 6.1$\pm$0.0e+01 \\
ELM-FBPINN & (10, [1, 10, 20]) & (1, 16) & 2.9 & 5.9e+16 & LSQR & 7.7$\pm$2.1e-03 & 1.5$\pm$0.0e+02 \\
ELM-FBPINN & (12, [1, 12, 24]) & (1, 16) & 2.9 & 6.6e+16 & LSQR & 1.3$\pm$0.4e-02 & 1.6$\pm$0.0e+02 \\
ELM-FBPINN & (14, [1, 14, 28]) & (1, 16) & 2.9 & 7.2e+16 & LSQR & 1.5$\pm$0.8e-02 & 2.0$\pm$0.0e+02 \\
ELM-FBPINN & (16, [1, 16, 32]) & (1, 16) & 2.9 & 8.0e+16 & LSQR & 7.8$\pm$1.6e-03 & 3.3$\pm$0.1e+02 \\
\bottomrule
\end{tabular}

 }
\captiontabscaling{$\omega$}{2D inhomogenous Helmholtz equation}{along each input dimension for}
\label{tab:ch3-helmholtz-scaling-results}
\end{table}

\rev{\paragraph{Interpretation}
The observed behaviour highlights the balance between problem complexity and model capacity. While increasing $J$ compensates for higher frequencies by maintaining local resolution, the global system becomes larger and more challenging to solve. The ELM-FBPINN remains more robust in this regime, suggesting that its formulation is better suited to handling the increased spectral content. The benefits of multilevel decomposition indicate that combining scales can further stabilise the approximation, particularly for highly oscillatory solutions.}

\rev{
\subsubsection{Accuracy--efficiency trade-off}
All 2D results highlight a clear trade-off between accuracy and computational efficiency for the ELM-FBPINN compared to the FBPINN. While the ELM-FBPINN initially converges significantly faster in terms of number of iteration steps and achieves nearly two orders of magnitude lower total runtime, the FBPINN attains lower final $e_{L^1}$ values in two-dimensional multi-scale settings. This behaviour reflects the underlying difference in model flexibility. In particular the FBPINN optimises all network parameters within each subdomain, providing greater expressivity and allowing finer adaptation to multi-scale structure, at the cost of slower and more computationally expensive optimisation. From a practical perspective, this suggests that ELM-FBPINN is particularly advantageous in regimes where fast convergence, robustness, and scalability are critical, such as large-scale or high-dimensional problems where gradient-based optimisation becomes challenging. Conversely, FBPINN may be preferred when the highest possible accuracy is required, especially in settings where fine-scale interactions dominate and additional computational cost is acceptable.
}

\subsubsection{Varying parameter initialisation}

Finally, we vary the weight scaling parameter $R$ in the LeCun initialisation of hidden basis function weights and biases of the baseline $L=1$ model \rev{for all problems tested} to ascertain the sensitivity of the ELM-FBPINN models to parameter initialisation.

\begin{figure}[t]
    \centering
    \includegraphics[width=0.9\linewidth]{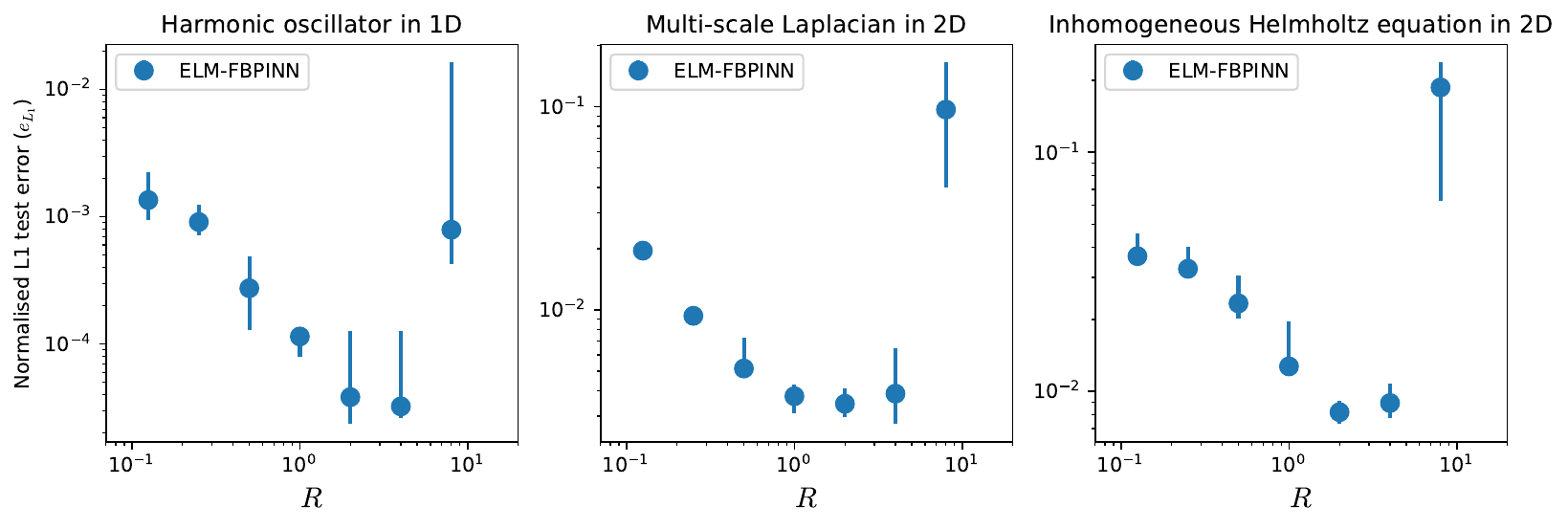}
    \caption{Final \rev{test error, measured in the $e_{L^1}$ norm,} when varying weight scaling values, $R$, \rev{for all problems tested}. \rev{The median final test error is denoted by a circle, and the range of final test errors by a vertical bar.}}
    \label{fig:weight-scaling}
\end{figure}

\paragraph{Key observations}
\Cref{fig:weight-scaling} shows that whilst the performance of ELM-FBPINN does depend on the choice of the weight scaling parameter, $R$, good performance is achieved for $R \in [1,5]$ for both the 1D and 2D problems. Indeed, even for $R$ outside this range, the \rev{test error} obtained is typically within one order of magnitude of the optimal value, indicating robustness with respect to the choice of parameter initialisation.

\section{Discussion and conclusions}
\label{sec:ch3_disc}

We compared PINN, FBPINN, ELM-FBPINN, and multilevel ELM-FBPINN across controlled ablation studies on a \rev{1D damped harmonic oscillator, 2D multi-scale Laplacian, and 2D inhomogeneous Helmholtz equation problem}. The experiments isolated the effects of localisation ($J$), width ($K$) and spectral complexity, while measuring convergence behaviour, \rev{error}, conditioning, and runtime.

\begin{table}[h]
\centering
\small
\renewcommand{\arraystretch}{1.2}
\begin{tabular}{lcccc}
\toprule
 & PINN & FBPINN & ELM-FBPINN & Multilevel ELM-FBPINN \\
\midrule
Multi-scale convergence & Poor & Stable & Stable (fast initial) & Stable (fast initial) \\
Final \rev{error} (1D) & \rev{Highest} & \rev{Low} & \rev{Lowest} & - \\
Final \rev{error} (2D) & \rev{Highest} & \rev{Lowest} & \rev{Low} & - \\
Scalability (frequency) & Poor & Stable & Stable & Stable \\
Runtime efficiency & Slow & Moderate & Fastest & Fast \\
\bottomrule
\end{tabular}
\caption{Qualitative comparison across experiments.}
\label{tab:summary-takeaways}
\end{table}
\vspace{-0.2cm}
\paragraph{Convergence behaviour}
Across all configurations, ELM-FBPINN exhibits a characteristic two-phase convergence: a rapid initial reduction of error followed by slower refinement. This reflects the structure of the global least-squares solve, where dominant residual components are eliminated early. Standard PINNs fail to converge on the multi-scale configurations considered, confirming that global gradient-based training alone is insufficient for oscillatory problems.

\paragraph{Accuracy trade-offs}
In 1D, ELM-FBPINN achieves the lowest errors in most ablations, except when activation saturation limits representational quality. In 2D, FBPINN generally attains lower final $e_{L^1}$ values, indicating that adaptive hidden weights are advantageous when resolving interacting spatial frequencies. Nevertheless, ELM-FBPINN consistently achieves errors within the same order of magnitude while requiring substantially less computational effort.

\paragraph{Scalability}
When frequency and model capacity are increased simultaneously, both FBPINN and ELM-FBPINN methods maintain stable error levels without catastrophic degradation. Furthermore, multilevel ELM-FBPINN perform best whilst only incurring slightly more computational cost than ELM-FBPINNs, suggesting that adding multiple coarser levels facilitates global subdomain communication and aids convergence.

\paragraph{Conditioning and redundancy}
The ELM-FBPINN least-squares matrix $\mathbf{M}$ typically exhibits large SVD-based condition numbers, often exceeding $10^{18}$ for large $K$. This growth is driven primarily by redundancy among randomly initialised basis functions. Increasing $K$ introduces near-linear dependencies, while excessive localisation ($J$ large) can also induce effective rank saturation when local solutions become simple. Despite this ill-conditioning, the least-squares solves remain sufficiently stable to produce accurate approximations, although redundancy limits further gains beyond moderate capacity.

\paragraph{Overall assessment}
Domain decomposition is essential for multi-scale problems. Within this framework, replacing gradient-based subdomain training by a linear least-squares solve yields substantial gains in convergence speed and competitive accuracy. While fully trainable hidden representations remain advantageous in higher-dimensional settings, ELM-FBPINN offers a favourable efficiency–accuracy trade-off and demonstrates robust scalability across the tested regimes, especially when multiple levels are employed. These results position random-feature-based domain decomposition as a viable and computationally efficient alternative to gradient-trained FBPINNs for structured scientific machine learning problems.

\section{Code and data availability}
\label{app:code_availability}
All code required to reproduce this work is available at \href{https://github.com/benmoseley/FBPINNs/tree/elm-paper/elm-paper}{https://github.com/benmoseley/FBPINNs/tree/elm-paper/elm-paper}. All data is synthetically generated using the provided code.

\section*{Funding declaration} 
The first author was supported by a University of Strathclyde Research Excellence Award (REA), funded jointly by the Department of Mathematics and Statistics, the Faculty of Science, the National Manufacturing Institute Scotland (NMIS), and the Advanced Forming Research Centre (AFRC).

\backmatter

\bibliography{paper_MLCSE}% common bib file

@inproceedings{yildiz2026fastmultiscale,
  title     = {Fast Multiscale {PDE} Solvers via Multilevel Domain Decomposition and Random Features},
  author    = {Yıldız, Eray and Datar, Chinmay and Rahma, Atamert and Dolean, Victorita and Dietrich, Felix},
  booktitle = {ICLR 2026 Workshop on AI for PDEs},
  year      = {2026},
  url       = {https://openreview.net/forum?id=wA6SsMmftX}
}

@inproceedings{datar2026frozenpinn,
  title     = {Frozen-PINN: Fast training of accurate physics-informed neural networks with temporal causality},
  author    = {Datar, Chinmay and Kapoor, Taniya and Chandra, Abhishek and Sun, Qing and Bolager, Erik Lien and Burak, Iryna and Veselovska, Anna and Fornasier, Massimo and Dietrich, Felix},
  booktitle = {International Conference on Learning Representations (ICLR)},
  year      = {2026},
  url       = {https://openreview.net/forum?id=3VdSuh3sie},
  note      = {ICLR 2026 Oral}
}

@article{vanBeek2026RRQR,
  title   = {Local feature filtering for scalable and well-conditioned domain-decomposed random feature methods},
  author  = {van Beek, Jan Willem and Dolean, Victorita and Moseley, Ben},
  journal = {Computer Methods in Applied Mechanics and Engineering},
  volume  = {449},
  pages   = {118583},
  year    = {2026},
  issn    = {0045-7825},
  doi     = {10.1016/j.cma.2025.118583},
  url     = {https://doi.org/10.1016/j.cma.2025.11858}
}

@article{Schiassi2021_ETFC,
title = {Extreme theory of functional connections: A fast physics-informed neural network method for solving ordinary and partial differential equations},
journal = {Neurocomputing},
volume = {457},
pages = {334-356},
year = {2021},
issn = {0925-2312},
doi = {https://doi.org/10.1016/j.neucom.2021.06.015},
url = {https://www.sciencedirect.com/science/article/pii/S0925231221009140},
author = {Enrico Schiassi and Roberto Furfaro and Carl Leake and Mario {De Florio} and Hunter Johnston and Daniele Mortari}
}

@article{DongLi2021_ELM-DDM,
title = {Local extreme learning machines and domain decomposition for solving linear and nonlinear partial differential equations},
journal = {Computer Methods in Applied Mechanics and Engineering},
volume = {387},
pages = {114129},
year = {2021},
issn = {0045-7825},
doi = {https://doi.org/10.1016/j.cma.2021.114129},
url = {https://www.sciencedirect.com/science/article/pii/S0045782521004606},
author = {Suchuan Dong and Zongwei Li}
}

@article{pao1992functional,
	title={Functional-link net computing: theory, system architecture, and functionalities},
	author={Pao, Y-H and Takefuji, Yoshiyasu},
	journal={Computer},
	volume={25},
	number={5},
	pages={76--79},
	year={1992},
	publisher={IEEE},
    doi={10.1109/2.144401}
}

@ARTICLE{2020SciPy-NMeth,
  author  = {Virtanen, Pauli and Gommers, Ralf and Oliphant, Travis E. and
            Haberland, Matt and Reddy, Tyler and Cournapeau, David and
            Burovski, Evgeni and Peterson, Pearu and Weckesser, Warren and
            Bright, Jonathan and {van der Walt}, St{\'e}fan J. and
            Brett, Matthew and Wilson, Joshua and Millman, K. Jarrod and
            Mayorov, Nikolay and Nelson, Andrew R. J. and Jones, Eric and
            Kern, Robert and Larson, Eric and Carey, C J and
            Polat, {\.I}lhan and Feng, Yu and Moore, Eric W. and
            {VanderPlas}, Jake and Laxalde, Denis and Perktold, Josef and
            Cimrman, Robert and Henriksen, Ian and Quintero, E. A. and
            Harris, Charles R. and Archibald, Anne M. and
            Ribeiro, Ant{\^o}nio H. and Pedregosa, Fabian and
            {van Mulbregt}, Paul and {Contributors}},
  title   = {{{SciPy} 1.0: Fundamental Algorithms for Scientific
            Computing in Python}},
  journal = {Nature Methods},
  year    = {2020},
  volume  = {17},
  pages   = {261--272},
  adsurl  = {https://rdcu.be/b08Wh},
  doi     = {10.1038/s41592-019-0686-2},
}

@article{Shang_Heinlein_Mishra_Wang_2025, title={Overlapping {S}chwarz preconditioners for randomized neural networks with domain decomposition}, volume={442}, DOI={10.1016/j.cma.2025.118011}, number={1}, journal={Computer Methods in Applied Mechanics and Engineering}, author={Shang, Yong and Heinlein, Alexander and Mishra, Siddhartha and Wang, Fei}, year={2025}, month=jul }

@article{Jaeger:2009:RCA,
title = {Reservoir computing approaches to recurrent neural network training},
journal = {Computer Science Review},
volume = {3},
number = {3},
pages = {127-149},
year = {2009},
issn = {1574-0137},
doi = {https://doi.org/10.1016/j.cosrev.2009.03.005},
url = {https://www.sciencedirect.com/science/article/pii/S1574013709000173},
author = {Mantas Lukoševičius and Herbert Jaeger}
}

@article{jagtap2020conservative,
	title={Conservative physics-informed neural networks on discrete domains for conservation laws: Applications to forward and inverse problems},
	author={Jagtap, Ameya D and Kharazmi, Ehsan and Karniadakis, George Em},
	journal={Computer Methods in Applied Mechanics and Engineering},
	volume={365},
	pages={113028},
	year={2020},
    doi={10.1016/j.cma.2020.113028},
	publisher={Elsevier}
}

@article{jagtap2020extended,
	title={Extended physics-informed neural networks ({XPINNs}): A generalized space-time domain decomposition based deep learning framework for nonlinear partial differential equations},
	author={Jagtap, Ameya D and Karniadakis, George Em},
	journal={Communications in Computational Physics},
	volume={28},
	number={5},
	year={2020},
    pages={2002--2041},
	publisher={Brown Univ., Providence, RI (United States)},
    doi={ 10.4208/cicp.OA-2020-0164}
}

@article{chen:2022:BTM,
	title={Bridging traditional and machine learning-based algorithms for solving {PDE}s: {T}he random feature method},
	author={Chen, Jingrun and Chi, Xurong and E, Weinan and Yang, Zhouwang},
	journal={Journal of Machine Learning},
	volume={1},
    number={3},
	pages={268--98},
	year={2022},
    doi={10.4208/jml.220726}
}

@article{Dolean:MDD:2024,
title = {Multilevel domain decomposition-based architectures for physics-informed neural networks},
journal = {Computer Methods in Applied Mechanics and Engineering},
volume = {429},
pages = {117116},
year = {2024},
issn = {0045-7825},
author = {V. Dolean and A. Heinlein and S. Mishra and Ben Moseley},
doi = {10.1016/j.cma.2024.117116}
}

@inproceedings{dolean2022finite,
	title={Finite basis physics-informed neural networks as a {Schwarz} domain decomposition method},
	author={Dolean, Victorita and Heinlein, Alexander and Mishra, Siddhartha and Moseley, Ben},
	booktitle={International Conference on Domain Decomposition Methods},
	pages={165--172},
	year={2022},
	organization={Springer},
    doi={10.1007/978-3-031-50769-4_19}
}

@misc{bradbury_jax_2018,
	title = {{JAX}: composable transformations of {Python}+{NumPy} programs},
	url = {http://github.com/google/jax},
	author = {Bradbury, James and Frostig, Roy and Hawkins, Peter and Johnson, Matthew James and Leary, Chris and Maclaurin, Dougal and Necula, George and Paszke, Adam and VanderPlas, Jake and Wanderman-Milne, Skye and Zhang, Qiao},
	year = {2018},
}

@article{Moseley2023,
author = {Moseley, Ben and Markham, Andrew and Nissen-Meyer, Tarje},
doi = {10.1007/S10444-023-10065-9},
isbn = {0123456789},
issn = {1572-9044},
journal = {Advances in Computational Mathematics 2023 49:4},
month = {jul},
number = {4},
pages = {1--39},
publisher = {Springer},
title = {{Finite basis physics-informed neural networks ({FBPINNs}): a scalable domain decomposition approach for solving differential equations}},
url = {https://link.springer.com/article/10.1007/s10444-023-10065-9},
volume = {49},
year = {2023}
}

@article{Dwivedi2021,
author = {Dwivedi, Vikas and Parashar, Nishant and Srinivasan, Balaji},
doi = {10.1016/j.neucom.2020.09.006},
issn = {18728286},
journal = {Neurocomputing},
month = {jan},
pages = {299--316},
publisher = {Elsevier B.V.},
title = {{Distributed learning machines for solving forward and inverse problems in partial differential equations}},
volume = {420},
year = {2021}
}

@article{Huang2006,
author = {Huang, Guang Bin and Zhu, Qin Yu and Siew, Chee Kheong},
doi = {10.1016/j.neucom.2005.12.126},
issn = {09252312},
journal = {Neurocomputing},
month = {dec},
number = {1-3},
pages = {489--501},
publisher = {Elsevier},
title = {{Extreme learning machine: Theory and applications}},
volume = {70},
year = {2006}
}

@inproceedings{Rahaman2018,
author = {Rahaman, Nasim and Baratin, Aristide and Arpit, Devansh and Draxler, Felix and Lin, Min and Hamprecht, Fred and Bengio, Yoshua and Courville, Aaron},
booktitle = {36th International Conference on Machine Learning, ICML 2019},
pages = {9230--9239},
title = {On the spectral bias of neural networks},
volume = {97},
year = {2019},
organization = {International Machine Learning Society (IMLS)},
url = {https://proceedings.mlr.press/v97/rahaman19a.html}
}

@article{Lagaris1998,
author = {Lagaris, Isaac Elias and Likas, Aristidis and Fotiadis, Dimitrios I.},
doi = {10.1109/72.712178},
eprint = {9705023},
issn = {10459227},
journal = {IEEE Transactions on Neural Networks},
number = {5},
pages = {987--1000},
title = {{Artificial neural networks for solving ordinary and partial differential equations}},
volume = {9},
year = {1998}
}

@article{Raissi2019,
author = {Raissi, M. and Perdikaris, P. and Karniadakis, G. E.},
doi = {10.1016/j.jcp.2018.10.045},
issn = {10902716},
journal = {Journal of Computational Physics},
pages = {686--707},
publisher = {Elsevier Inc.},
title = {{Physics-informed neural networks: A deep learning framework for solving forward and inverse problems involving nonlinear partial differential equations}},
url = {https://doi.org/10.1016/j.jcp.2018.10.045},
volume = {378},
year = {2019}
}

@article{virtanen_scipy_2020,
	title = {{SciPy} 1.0: fundamental algorithms for scientific computing in {Python}},
	volume = {17},
	copyright = {2020 The Author(s)},
	issn = {1548-7105},
	shorttitle = {{SciPy} 1.0},
	url = {https://www.nature.com/articles/s41592-019-0686-2},
	doi = {10.1038/s41592-019-0686-2},
	language = {en},
	number = {3},
	urldate = {2026-01-19},
	journal = {Nature Methods},
	publisher = {Nature Publishing Group},
	author = {Virtanen, Pauli and Gommers, Ralf and Oliphant, Travis E. and Haberland, Matt and Reddy, Tyler and Cournapeau, David and Burovski, Evgeni and Peterson, Pearu and Weckesser, Warren and Bright, Jonathan and van der Walt, St{\textbackslash}'\{e\}fan J. and Brett, Matthew and Wilson, Joshua and Millman, K. Jarrod and Mayorov, Nikolay and Nelson, Andrew R. J. and Jones, Eric and Kern, Robert and Larson, Eric and Carey, C. J. and Polat, Ilhan and Feng, Yu and Moore, Eric W. and VanderPlas, Jake and Laxalde, Denis and Perktold, Josef and Cimrman, Robert and Henriksen, Ian and Quintero, E. A. and Harris, Charles R. and Archibald, Anne M. and Ribeiro, Ant{\textbackslash}{\textasciicircum}\{o\}nio H. and Pedregosa, Fabian and van Mulbregt, Paul},
	month = mar,
	year = {2020},
	pages = {261--272},
}

@inproceedings{lecun_efficient_1998,
author="LeCun, Yann A.
and Bottou, L{\'e}on
and Orr, Genevieve B.
and M{\"u}ller, Klaus-Robert",
editor="Montavon, Gr{\'e}goire
and Orr, Genevi{\`e}ve B.
and M{\"u}ller, Klaus-Robert",
title="Efficient {BackProp}",
booktitle="Neural Networks: Tricks of the Trade",
edition="Second",
year="2012",
publisher="Springer Berlin Heidelberg",
address="Berlin, Heidelberg",
pages="9--48",
doi="10.1007/978-3-642-35289-8_3"
}

@misc{Anderson:2024:ELM,
	title = {{ELM}-{FBPINN}: efficient finite-basis physics-informed neural networks},
	shorttitle = {{ELM}-{FBPINN}},
	url = {http://arxiv.org/abs/2409.01949},
	author = {Anderson, Samuel and Dolean, Victorita and Moseley, Ben and Pestana, Jennifer},
	month = sep,
	year = {2024},
}

\end{document}